\documentclass[a4paper, final, 10pt]{article}

\usepackage{lipsum}
\usepackage{amsmath, amssymb, latexsym, amscd,amsfonts,amstext}
\usepackage[mathscr]{eucal}
\usepackage{graphicx}
\usepackage{subfig}
\usepackage{float}
\usepackage{color}
\usepackage{hyperref}
\usepackage[utf8]{inputenc}
\usepackage[english]{babel}
\usepackage{pgfplotstable}
\usepackage{tikz}
\usepackage{epstopdf}
\usepackage{algorithmic}
\usepackage{amsopn}

\newtheorem{theorem}{Theorem}

\numberwithin{theorem}{section}
\numberwithin{equation}{section}

 \ifpdf
  \DeclareGraphicsExtensions{.eps,.pdf,.png,.jpg}
\else
  \DeclareGraphicsExtensions{.eps}
\fi 

\begin{document}


%
%
%
%

\title{Convexification of a 3-D coefficient inverse scattering problem\thanks{Supported by US Army Research Laboratory and US Army Research Office grant W911NF-15-1-0233 and by the Office of Naval Research grant N00014-15-1-2330. In addition, the work of Kolesov A.E. was partially supported by Mega-grant of the Russian Federation Government (N14.Y26.31.0013) and RFBR (project N17-01-00689A)}}
\author{Michael V. Klibanov \thanks{The corresponding author} \thanks{Department of Mathematics \& Statistics, University of North Carolina at Charlotte, Charlotte, NC 28223, USA (mklibanv@uncc.edu, 
akolesov@uncc.edu)},  
Aleksandr E.\ Kolesov \footnotemark[3] \thanks{Institute of Mathematics and Information Science, North-Eastern Federal University, Yakutsk, Russia (ae.kolesov@s-vfu.ru)}}

\maketitle

\begin{abstract}
A version of the so-called \textquotedblleft convexification" numerical
method for a coefficient inverse scattering problem for the 3D Hemholtz
equation is developed analytically and tested numerically. Backscattering
data are used, which result from a single direction of the propagation of
the incident plane wave on an interval of frequencies. The method converges
globally. The idea is to construct a weighted Tikhonov-like functional. The
key element of this functional is the presence of the so-called Carleman
Weight Function (CWF). This is the function which is involved in the Carleman
estimate for the Laplace operator. This functional is strictly convex on any
appropriate ball in a Hilbert space for an appropriate choice of the
parameters of the CWF. Thus, both the absence of local minima and
convergence of minimizers to the exact solution are guaranteed. Numerical tests demonstrate
a good performance of the resulting algorithm. Unlikeprevious the 
so-called tail functions globally convergent method, we neither do not 
impose the smallness assumption of the interval of wavenumbers, nor
we do not iterate with respect to the so-called tail functions.
\end{abstract}




\section{Introduction}

\label{sec:1}

In this work, we develop a version of the so-called \textquotedblleft
convexification" numerical method for a coefficient inverse scattering
problem (CISP) for the 3D Helmholtz equation with backscattering data
resulting from a single measurement event which is generated by a single
direction of the propagation of the incident plane wave on an interval of
frequencies. We present both the theory and numerical results. Our method
converges globally. This is a generalization to the 3D case of our (with
coauthors) previous 1D version of the convexification \cite%
{KlibanovKolesov17}. Three main advantages of the convexification method
over the previously developed the so-called \textquotedblleft tail
functions" globally convergent method for a similar CISP \cite%
{BeilinaKlibanov08,
BeilinaKlibanov12,KlibanovLiem16,KlibanovKolesov17exp,KlibanovLiem17buried,KlibanovLiem17exp}
are: (1) To solve our problem, we construct a globally strictly convex cost
functional with the Carleman Weight Function (CWF) in it, (2) we do not
impose in our convergence analysis the smallness assumption on the interval
of wavenumbers, and (3) we do not iterate with respect to the so-called
\textquotedblleft tail functions".

It is well known that any CISP is both highly nonlinear and ill-posed. These
two factors cause substantial difficulties in numerical solutions of these
problems. A \emph{globally convergent method }(GCM) for a CISP is such a
numerical method, which has a rigorous guarantee of reaching a sufficiently
small neighborhood of the exact solution of that CISP without any advanced
knowledge of this neighborhood. In addition, the size of this neighborhood
should depend only on approximation errors and the level of noise in the
data.

Over the years the first author with coauthors has proposed a variety of
globally convergent methods for CISPs with single measurement data, see,
e.g. \cite{BeilinaKlibanov08, BeilinaKlibanov12, KlibanovLiem16,
KlibanovKolesov17exp, KlibanovLiem17exp,
BeilinaKlibanov15,KlibanovIoussoupova95, Klibanov97a, Klibanov97b,
KlibanovTimonov04, KlibanovKamburg16, KlibanovLoc16, KlibanovThanh15}, and
references cited therein. These methods can be classified into two types.
Methods of the first type, which we call the \emph{tail functions} methods,
are certain iterative processes. On each iterative step one solves the
Dirichlet boundary value problem for a linear elliptic Partial Differential
Equation (PDE). This PDE depends on the iteration number. The solution of
that problem enables one to update the unknown coefficient. Using this
update, one updates the so-called tail function, which is a complement of a
certain truncated integral, where the integration is carried out with
respect to the wavenumber. The stopping criterion for the iterative process
is developed computationally. The tail function method was successfully
tested on experimental backscattering data The tail function method was
successfully tested on experimental backscattering data \cite{
KlibanovLiem16, KlibanovKolesov17exp, KlibanovLiem17buried,
KlibanovLiem17exp, KlibanovLoc16}.

Globally convergent numerical methods of the second type are called the 
\emph{convexification} methods. They are based on the minimization of the
weighted Tikhonov-like functional with the CWF in it. The CWF is the
function which is involved in the Carleman estimate for the corresponding
PDE operator. The CWF can be chosen in such a way that the above functional
becomes strictly convex on a ball of an arbitrary radius in a certain
Hilbert space (see some details in this section below). Note that the
majority of known numerical methods of solutions of nonlinear ill-posed
problems minimize conventional least squares cost functionals \cite%
{Chavent09,Goncharsky13,Goncharsky17}, which are usually non convex and have
multiple local minima and ravines, see, e.g. \cite{Scales92} for a good
numerical example of multiple local minima. Hence, a gradient-like method
for such a functional converges to the exact solution only if the starting
point of iterations is located in a sufficiently small neighborhood of this
solution. Some other effective approache to numerical methods for nonlinear
ill-posed problems can be found in \cite{Lakhal2,Lakhal3}.

Various versions of the convexification methods have been proposed since the
first work \cite{KlibanovIoussoupova95}, see \cite{Klibanov97a,Klibanov97b,
KlibanovTimonov04}. However, these versions have some theoretical gaps,
which have limited their numerical studies so far. In the recent works \cite%
{BeilinaKlibanov15, KlibanovKamburg16, KlibanovKoshev16} the attention to
the convexification method was revived. Theoretical gaps were eliminated in 
\cite{BakushinskiiKlibanov17} and thorough numerical studies for one
dimensional problems were performed \cite{KlibanovKolesov17,KlibanovThanh15}%
. Besides, in \cite{Klibanov15} the convexification method was developed for
ill-posed problems for quasilinear PDEs and corresponding numerical studies
for the 1D case were conducted in \cite%
{KlibanovKoshev16,BakushinskiiKlibanov17}. The idea of any version of the
convexification has direct roots in the method of \cite{KlibanovBukhgeim81},
which is based on Carleman estimates. The method of \cite{KlibanovBukhgeim81}
was originally designed only for proofs of uniqueness theorems for CIPs,
also see, e.g. the book \cite{KlibanovTimonov04} and the recent survey \cite%
{Ksurvey}. Recently an interesting version of the convexification was
published in \cite{Baud} for a CISP for the hyperbolic equation $%
u_{tt}=\Delta u+a\left( x\right) u$ with the unknown coefficient $a\left(
x\right) $ in the case when one of initial conditions does not non-vanish.
The method of \cite{Baud} is also based on the idea of \cite%
{KlibanovBukhgeim81} and has some roots in \cite%
{BeilinaKlibanov15,KlibanovKamburg16}.

By the convexification, one constructs a weighted Tikhonov-like functional $%
J_{\lambda }$ on a closed ball $\overline{B\left( R\right) }$ of an
arbitrary radius $R>0$ and with the center at $\left\{ 0\right\} $ in an
appropriate Hilbert space. Here $\lambda >0$ is a parameter. The key theorem
claims that one can choose a number $\lambda \left( R\right) >0$ such that
for all $\lambda \geq \lambda \left( R\right) $ the functional $J_{\lambda }$
is strictly convex on $\overline{B\left( R\right) }.$ Furthermore, the
existence of the unique minimizer of $J_{\lambda }$ on $\overline{B\left(
R\right) }$ as well as convergence of minimizers to the exact solution when
the level of noise in the data tends to zero are proven. In addition, it is
proven that the gradient projection method reaches a sufficiently small
neighborhood of the exact coefficient when starting from an arbitrary point
of $B\left( R\right) $. Since $R>0$ is an arbitrary number, then this is a 
\emph{globally convergent} numerical method.

Due to a broad variety of applications, Inverse Scattering Problems (ISPs)
are quite popular in the community of experts in inverse problems. There are
plenty of works dedicated to this topic. Since this paper is not a survey,
we refer to only few of them, e.g. \cite%
{Goncharsky13,Goncharsky17,Lakhal2,Lakhal3,Am1,Am2,Am3,Bao,Buhan,Chow1,
Chow2,Ito,Jin,Kab1,Kab,Kab3,Lakhal1,Liu1,Liu2,Liu3} and references cited
thereein. We note that the authors of \cite{Chow1} have considered a
modified tail functions method. As stated above, we are interested in a CISP
for the Helmholtz equation with the data generated by a single measurement
event. As to the CISPs with multiple measurements, we refer to a global
reconstruction procedure, which was developed and numerically implemented in 
\cite{Kab1}, also see \cite{Kab,Kab3} for further developments and numerical
studies. Actually, this is an effective extension of the classical 1D
Gelfand-Krein-Levitan method on the 2D case.

In section 2 we formulate our forward and inverse problems. In section 3 we
construct the weighted Tikhonov-like functional with the CWF in it. In
section 4 we formulate our theorems. \ We prove them in section 5. In
section 6 we present numerical results.

\section{Problem Statement}

\label{sec:2}

\subsection{The Helmholtz equation}

\label{sec:2.1}

Just as in the majority of the above cited previous works of the first
author with coauthors about GCM, we focus in this paper applications to the
detection and identification of targets, which mimic antipersonnel land
mines (especially plastic mines, i.e. dielectrics) and improvised explosive
devices (IEDs) using measurements of a single component of the electric wave
field. In this case the medium is assumed to be non magnetic, non absorbing,
and the dielectric constant in it should be represented by a function, which
is mostly a constant with some small sharp inclusions inside (however, we do
not assume in our theory such a structure of the dielectric constant). These
inclusions model antipersonnel land mines and IEDs. Suppose that the
incident electric field has only one non zero component. It was established
numerically in \cite{Liem} that the propagation of that component through
such a medium is well governed by the Helmholtz equation rather than by the
full Maxwell's system. Besides, in all above cited works of the first author
with coauthors about experimental data those targets were accurately imaged
by the above mentioned tail functions GCM using experimentally measured
single component of the electric field and modeling the propagation of that
component by the Helmholtz equation. In addition, we are unaware about a GCM
for a CISP with single measurement data for the Maxwell's system. Thus, we
use the Helmholtz equation below.

The need of the detection and identification of, e.g. land mines, might, in
particular, occur on a battlefield. Due to the security considerations, the
amount of collected data should be small in this case, and these should be
the backcattering data. Thus, we use only a single direction of the
propagation of the incident plane wave of the electric field and assume
measurements of only the backscattering part of the corresponding component
of that field.

\subsection{Forward and inverse problems}

\label{sec:2.2}

Let $\mathbf{x}=(x,y,z)\in \mathbb{R}^{3}$. Let $b,d,\xi >0$ be three
numbers. It is convenient for our numerical studies (section 6) to define
from the beginning the domain of interest $\Omega $ and the backscattering
part $\Gamma $ of its boundary as%
\begin{equation}
\Omega =\left\{ \left( x,y,z\right) :\left\vert x\right\vert ,\left\vert
y\right\vert <b,z\in \left( -\xi ,d\right) \right\} ,\text{ }\Gamma =\left\{
\left( x,y,z\right) :\left\vert x\right\vert ,\left\vert y\right\vert
<b,z=-\xi \right\} .  \label{eq:2.1}
\end{equation}%
Let the function $c(\mathbf{x})$ be the spatially distributed dielectric
constant and $k$ be the wavenumber. We consider the following forward
problem for the Helmholtz equation: 
\begin{equation}
\Delta u+k^{2}\,c(\mathbf{x})\,u=0,\quad \mathbf{x}\in \mathbb{R}^{3},
\label{eq:helmholtz}
\end{equation}%
\begin{equation}
u\left( \mathbf{x},k\right) =u_{s}\left( \mathbf{x},k\right) +u_{i}\left( %
\mathbb{x},k\right) ,  \label{eq:utotal}
\end{equation}%
where $u(\mathbf{x},k)$ is the total wave, $u_{s}(\mathbf{x},k)$ is the
scattered wave, and $u_{i}(\mathbf{x},k)$ is the incident plane wave
propagating along the positive direction of the $z-$axis, 
\begin{equation}
u_{i}(\mathbf{x},k)=e^{ikz}.  \label{eq:uinc}
\end{equation}%
The scattered wave $u_{s}(\mathbf{x},k)$ satisfies the Sommerfeld radiation
condition: 
\begin{equation}
\lim_{r\rightarrow \infty }r\left( \frac{\partial u_{s}}{\partial r}%
-iku_{s}\right) =0,\quad r=\left\vert \mathbf{x}\right\vert .
\label{eq:sommerfled}
\end{equation}%
Also, the function $c(\mathbf{x})$ satisfies the following conditions: 
\begin{equation}
c(\mathbf{x})=1+\beta \left( \mathbf{x}\right) ,\text{ }\beta \left( \mathbf{%
x}\right) \geq 0,\,\mathbf{x}\in \mathbb{R}^{3},\quad \mbox{and}c(\mathbf{x}%
)=1,\,\mathbf{x}\notin \overline{\Omega }.  \label{eq:coef}
\end{equation}%
The assumption of (\ref{eq:coef}) $c(\mathbf{x})=1$ in $\mathbb{R}%
^{3}\setminus \Omega $ means that we have vacuum outside of the domain $%
\Omega .$ Finally, we assume that $c(\mathbf{x})\in C^{15}(\mathbb{R}^{3})$.
This smoothness condition was imposed to derive the asymptotic behavior of
the solution of the Helmholtz equation (\ref{eq:helmholtz}) at $k\rightarrow
\infty $ \cite{KlibanovRomanov16}. We also note that extra smoothness
conditions are usually not of a significant concern when a CIP is
considered, see, e.g. theorem \ref{thm:4.1} in \cite{Rom}. In particular, this
smoothness condition implies that the function $u\left( \mathbf{x},k\right)
\in C^{16+\gamma }\left( \overline{G}\right) ,\forall \gamma \in \left(
0,1\right) ,\forall k>0,$ where $C^{16+\gamma }\left( \overline{G}\right) $
is the H\"{o}lder space and $G\subset \mathbb{R}^{3}$ is an arbitrary
bounded domain \cite{Gilbarg}. Also, it follows from lemma 3.3 of \cite%
{KlibanovLiem16} that the derivative $\partial _{k}u\left( \mathbf{x}%
,k\right) $ exists for all $\mathbf{x}\in \mathbb{R}^{3},k>0$ and satisfies
the same smoothness condition as the function $u\left( \mathbf{x},k\right) .$

\textbf{Coefficient Inverse Scattering Problem (CISP).} \emph{Let the domain 
$\Omega $ and the backscattering part $\Gamma \subset \partial \Omega $ of
its boundary be as in (\ref{eq:2.1}). Let the wavenumber $k\in \lbrack 
\underline{k},\overline{k}],$ where }$[\underline{k},\overline{k}]$\emph{$%
\subset \left( 0,\infty \right) $ is an interval of wavenumbers. Determine
the function $c(\mathbf{x}),\,\mathbf{x}\in \Omega $, assuming that the
following function $g(\mathbf{x},k)$ is given}: 
\begin{equation}
u(\mathbf{x},k)=g_{0}(\mathbf{x},k),\quad \mathbf{x}\in \Gamma ,\,k\in
\lbrack \underline{k},\overline{k}].  \label{eq:cisp}
\end{equation}

In addition to the data (\ref{eq:cisp}) we can obtain the boundary
conditions for the derivative of the function $u(\mathbf{x},k)$ in the $z-$%
direction using the data propagation procedure (section 6.2), 
\begin{equation}
u_{z}(\mathbf{x},k)=g_{1}(\mathbf{x},k),\quad \mathbf{x}\in \Gamma ,\,k\in
\lbrack \underline{k},\overline{k}].  \label{eq:gz0}
\end{equation}

In addition, we complement Dirichlet (\ref{eq:cisp}) and Neumann (\ref%
{eq:gz0}) boundary conditions on $\Gamma $ with the heuristic Dirichlet
boundary condition at the rest of the boundary $\partial \Omega $ as:%
\begin{equation}
u(\mathbf{x},k)=e^{ikz},\mathbf{x}\in \partial \Omega \diagdown \Gamma
,\,k\in \lbrack \underline{k},\overline{k}].  \label{2.2}
\end{equation}%
This boundary condition coincides with the one for the uniform medium with $%
c\left( \mathbf{x}\right) \equiv 1.$ To justify (\ref{2.2}), we recall that,
using the tail functions method, it was demonstrated in sections 7.6 and 7.7
of \cite{KlibanovLiem16} that (\ref{2.2}) does not affect much the
reconstruction accuracy as compared with the correct Dirichlet boundary
condition. Besides, (\ref{2.2}) has always been used in works \cite%
{KlibanovLiem16, KlibanovKolesov17exp, KlibanovLiem17buried,
KlibanovLiem17exp} with experimental data, where accurate results were
obtained by the tail functions GCM.

The uniqueness of the solution of this CISP is an open and long standing
problem. In fact, uniqueness of a similar coefficient inverse problem can be
currently proven only in the case if the right hand side of equation (\ref%
{eq:helmholtz}) is a function which is not vanishing in $\overline{\Omega }.$
This can be done by the method of \cite%
{KlibanovBukhgeim81,KlibanovTimonov04,Ksurvey}. Hence, we assume below the
uniqueness of our CISP.

\subsection{Travel time}

\label{sec:2.3}

The Riemannian metric generated by the function $c(\mathbf{x})$ is: 
\begin{equation*}
d\tau (\mathbf{x})=\sqrt{c(\mathbf{x})}|d\mathbf{x}|,\quad |d\mathbf{x}|=%
\sqrt{(dx)^{2}+(dy)^{2}+(dz)^{2}}.
\end{equation*}%
Fix the number $a>0.$ Consider the plane $P_{a}=\{(x,y,-a):x,y\in \mathbb{R}%
\}.$ We assume that $\Omega \subset \left\{ z>-a\right\} $ and impose
everywhere below the following condition on the function $c(\mathbf{x})$:

\textbf{Regularity Assumption}. \emph{For any point }$x\in \mathbb{R}^{3}$%
\emph{\ there exists a unique geodesic line }$\Gamma (x,a)$\emph{, with
respect to the metric }$d\tau $\emph{, connecting }$x$\emph{\ with the plane 
}$P_{a}$\emph{\ and perpendicular to }$P_{a}$\emph{.}

A sufficient condition of the regularity of geodesic lines is \cite{Rom3}: 
\begin{equation*}
\sum_{i,j=1}^{3}\frac{\partial ^{2}c\left( \mathbf{x}\right) }{%
\partial x_{i}\partial x_{j}}\xi _{i}\xi _{j}\geq 0,\forall \mathbf{x}\in 
\overline{\Omega },\forall \mathbf{\xi }\in \mathbb{R}^{3}.
\end{equation*}

We introduce the travel time $\tau (\mathbf{x})$ from the plane $P_{a}$ to
the point $\mathbf{x}$ as \cite{KlibanovRomanov16} 
\begin{equation*}
\tau (\mathbf{x})=\int_{\Gamma (\mathbf{x},a)}\sqrt{c\left( \mathbf{%
\xi }\right) }d\sigma .
\end{equation*}

\section{The Weighted Tikhonov Functionals}

\label{sec:3}

\subsection{The asymptotic behavior}

\label{sec:3.1}

It was proven in \cite{KlibanovRomanov16} that the following asymptotic
behavior of the function $u(\mathbf{x},k)$ is valid: 
\begin{equation}
u(\mathbf{x},k)=A(\mathbf{x})e^{ik\tau (\mathbf{x})}\left[ 1+s\left( \mathbf{%
x},k\right) \right] ,\text{ }\mathbf{x}\in \overline{\Omega },k\rightarrow
\infty ,  \label{eq:uasymptotics}
\end{equation}%
where the function $s\left( \mathbf{x},k\right) $ is such that 
\begin{equation}
s\left( \mathbf{x,}k\right) =O\left( \frac{1}{k}\right) ,\partial
_{k}s\left( \mathbf{x,}k\right) =O\left( \frac{1}{k}\right) ,\text{ }\mathbf{%
x}\in \overline{\Omega },k\rightarrow \infty .  \label{3.1}
\end{equation}%
Here the function $A(\mathbf{x})>0$ and $\tau (\mathbf{x})$ is the length of
the geodesic line in the Riemannian metric generated by the function $c(%
\mathbf{x})$. Denote 
\begin{equation}
w(\mathbf{x},k)=\frac{u(\mathbf{x},k)}{u_{i}(\mathbf{x},k)}.  \label{eq:w}
\end{equation}

Using (\ref{eq:uasymptotics}), (\ref{3.1}) and (\ref{eq:w}), we obtain for $%
\mathbf{x}\in \overline{\Omega },k\rightarrow \infty $ that 
\begin{equation}
w(\mathbf{x},k)=A(\mathbf{x})e^{ik(\tau (\mathbf{x})-z)}\left[ 1+s\left( 
\mathbf{x},k\right) \right] .  \label{eq:wasymptotics}
\end{equation}%
Using (\ref{eq:uasymptotics}) and (\ref{eq:wasymptotics}), we uniquely
define the function $\log w(\mathbf{x},k)$ for $\mathbf{x}\in \Omega $, $%
k\in \lbrack \underline{k},\overline{k}]$ for sufficiently large values of $%
\underline{k}$ as 
\begin{equation}
\log w(\mathbf{x},k)=\ln A(\mathbf{x})+ik(\tau (\mathbf{x})-z)+%
\mathop{\displaystyle \sum }_{n=1}^{\infty }\frac{\left( -1\right)
^{n-1}}{n}\left( s(\mathbf{x} ,k)\right) ^{n}.  \label{3.2}
\end{equation}%
Obviously for so defined function $\log w(\mathbf{x},k)$ we have that $\exp %
\left[ \log w(\mathbf{x},k)\right] $ equals to the right hand side of (\ref%
{eq:wasymptotics}). Thus, we assume below that the number $\underline{k}$ is
sufficiently large.

\subsection{The integro-differential equation}

\label{sec:3.2}

It follows from (\ref{eq:helmholtz}), (\ref{eq:uinc}), (\ref{eq:coef}) and (%
\ref{eq:w}) that the function $w\left( \mathbf{x},k\right) $ satisfies the
following equation in the domain $\Omega $%
\begin{equation}
\Delta w+k^{2}\beta w+2ikw_{z}=0.  \label{eq:intdiffw}
\end{equation}

For $\mathbf{x}\in \Omega ,\,k\in \lbrack \underline{k},\overline{k}]$ we
define the function $v(\mathbf{x},k),$ 
\begin{equation}
v(\mathbf{x},k)=\frac{\log w(\mathbf{x},k)}{k^{2}}.  \label{eq:v}
\end{equation}%
Then%
\begin{equation}
\Delta v+k^{2}\left( \nabla v\right) ^{2}+2ikv_{z}+\beta (\mathbf{x})=0.
\label{eq:intdiffv}
\end{equation}%
Let $q(\mathbf{x},k)$ be the derivative of the function $v$ with respect to $%
k,$ 
\begin{equation}
q(\mathbf{x},k)=\partial _{k}v(\mathbf{x},k).  \label{eq:q}
\end{equation}%
Then 
\begin{equation}
v(\mathbf{x},k)=-\int_{k}^{\overline{k}}q\left( \mathbf{x},\kappa
\right) d\kappa +V(\mathbf{x}).  \label{eq:vq}
\end{equation}%
We call $V(\mathbf{x})$ the tail function: 
\begin{equation}
V(\mathbf{x})=v\left( \mathbf{x},\overline{k}\right) .  \label{3.30}
\end{equation}%
To eliminate the function $\beta (\mathbf{x})$ from equation (\ref%
{eq:intdiffv}), we differentiate (\ref{eq:intdiffv}) with respect to $k,$ 
\begin{equation}
\Delta q+2k\nabla v\cdot \left( k\nabla q+\nabla v\right) +2i\left(
kq_{z}+v_{z}\right) =0.  \label{eq:intdiffq}
\end{equation}%
Substituting (\ref{eq:vq}) into (\ref{eq:intdiffq}) leads to the following
integro-differential equationeq:

\begin{equation}
\begin{gathered} L(q)=\Delta q+2k\left( \nabla
V-\int_{k}^{\overline{k}}\nabla q(\mathbf{x},\kappa )d\kappa \right) \cdot
\left( k\nabla q+\nabla V-\int_{k}^{\overline{k}}\nabla q\left(
\mathbf{x},\kappa \right) d\kappa \right) \\ +2i\left(
kq_{z}+V_{z}-\int_{k}^{\overline{k}}q_{z}\left( \mathbf{x},\kappa \right)
d\kappa \right) =0. \end{gathered}  \label{eq:intdiff}
\end{equation}%
Finally, we complement this equation with the overdetermined boundary
conditions: 
\begin{equation}
\begin{gathered} q(\mathbf{x},k)=\phi _{0}(\mathbf{x},k),\quad
q_{z}(\mathbf{x},k) =\phi _{1}(\mathbf{x},k),\quad \mathbf{x}\in \Gamma
,\,k\in \lbrack \underline{k},\overline{k}], \\ q(\mathbf{x},k) =0,\quad
\mathbf{x}\in \partial \Omega \setminus \Gamma ,\,k\in \lbrack
\underline{k},\overline{k}], \end{gathered}  \label{eq:intdiffbcs}
\end{equation}%
where the functions $\phi _{0}$ and $\phi _{1}$ are calculated from the
functions $g_{0}$ and $g_{1}$ in (\ref{eq:cisp}), (\ref{eq:gz0}). The third
boundary condition (\ref{eq:intdiffbcs}) follows from (\ref{eq:uinc}), (\ref%
{2.2}), (\ref{eq:w}), (\ref{eq:v}) and (\ref{eq:q}).

Note that in (\ref{eq:intdiff}) both functions $q(\mathbf{x},k)$ and $V(%
\mathbf{x})$ are unknown. Hence, we approximate the function $V(\mathbf{x})$
first. Next, we solve the problem (\ref{eq:intdiff} ), (\ref{eq:intdiffbcs})
for the function $q(\mathbf{x},k)$.

\textbf{Remark 3.1.} Suppose that certain approximations for the functions $%
q(\mathbf{x},k)$\ and $V(\mathbf{x})$\ are found. Then an approximation for
the unknown coefficient $c\left( \mathbf{x}\right) $\ can be found via
backwards calculations: first, approximate the function $v\left( \mathbf{x}%
,k\right) $\ via (\ref{eq:vq}) and then approximate the function $\beta
\left( \mathbf{x}\right) $\ using equation (\ref{eq:intdiffv}) for a certain
value of $k\in \left[ \underline{k},\overline{k}\right] $. In our
computations we use $k=\underline{k}$ \ for that value of $k$. Next, one
should use (\ref{eq:coef}). Therefore, we focus below on approximating
functions $q(\mathbf{x},k)$ and $V(\mathbf{x}).$

\subsection{Approximation of the tail function}

\label{sec:3.3}

The method of this paper to approximate the tail function is different from
the method explored before in \cite{KlibanovKolesov17}. Also, unlike the
tail functions method, we do not update tails here.

It follows from (\ref{3.2}) and (\ref{3.30}) that there exists a function $p(%
\mathbf{x})$ such that 
\begin{equation}
v\left( \mathbf{x},k\right) =\frac{p\left( \mathbf{x}\right) }{k}+O\left( 
\frac{1}{k^{2}}\right) ,\quad q\left( \mathbf{x},k\right) =-\frac{p\left( 
\mathbf{x}\right) }{k^{2}}+O\left( \frac{1}{k^{3}}\right) ,\quad
k\rightarrow \infty ,\,\mathbf{x}\in \Omega .  \label{eq:vqasympt}
\end{equation}%
Since the number $\overline{k}$ is sufficiently large, we drop terms $%
O\left( 1/\overline{k}^{2}\right) $ and $O\left( 1/\overline{k}^{3}\right) $
in (\ref{eq:vqasympt}). Next, we approximately set 
\begin{equation}
v\left( \mathbf{x},k\right) =\frac{p\left( \mathbf{x}\right) }{k},\quad
q\left( \mathbf{x},k\right) =-\frac{p\left( \mathbf{x}\right) }{k^{2}},\quad
k\geq \overline{k},\,\mathbf{x}\in \Omega .  \label{eq:vqtail}
\end{equation}%
Substituting (\ref{eq:vqtail}) in (\ref{eq:intdiff}) and letting $k=%
\overline{k}$, we obtain 
\begin{equation}
\Delta V(\mathbf{x})=0,\quad \mathbf{x}\in \Omega .  \label{eq:tail}
\end{equation}%
This equation is supplemented by the following boundary conditions: 
\begin{equation}
V(\mathbf{x})=\psi _{0}(\mathbf{x}),\quad V_{z}(\mathbf{x})=\psi _{1}(%
\mathbf{x}),\quad \mathbf{x}\in \Gamma ,\quad V(\mathbf{x})=0,\quad \mathbf{x%
}\in \partial \Omega \setminus \Gamma ,  \label{eq:tailbcs}
\end{equation}%
where functions $\psi _{0}$ and $\psi _{1}$ can be computed using (\ref%
{eq:cisp}) and (\ref{eq:gz0}). Boundary conditions (\ref{eq:tailbcs}) are
over-determined ones. Due to the approximate nature of (\ref{eq:vqtail}), we
have observed that the obvious approach of finding the function $V(\mathbf{x}%
)$ by dropping the second boundary condition (\ref{eq:tailbcs}) and solving
the resulting Dirichlet boundary value problem for Laplace equation (\ref%
{eq:tail}) with the resulting boundary data (\ref{eq:tailbcs}) does not
provide satisfactory results. The same observation was made in \cite%
{KlibanovKolesov17} for the 1D case. Thus, we use a different approach to
approximate the function $V\left( \mathbf{x}\right) $.

Let the number $s>0$ be such that $s>\xi .$ Let $\lambda ,\nu >0$ be two
parameters which we will choose later. We introduce the CWF as%
\begin{equation}
\varphi _{\lambda }\left( z\right) =\exp \left[ 2\lambda \left( z+s\right)
^{-\nu }\right] ,  \label{3.3}
\end{equation}%
see Theorem \ref{thm:4.1} in section 4.1. Below we fix a number $\nu $ and allow $%
\lambda $ to change. We find an approximate solution of the problem (\ref%
{eq:tail}), (\ref{eq:tailbcs}) by minimizing the following cost functional
with the CWF in it:%
\begin{equation}
I_{\mu ,\alpha }\left( V\right) =\exp \left( -2\mu \left( s+d\right) ^{-\nu
}\right) \int_{\Omega }\left\vert \Delta V\right\vert ^{2}\varphi
_{\mu }\left( z\right) d\mathbf{x}+\alpha \Vert V\Vert _{H^{3}(\Omega )}^{2}.
\label{eq:Jtail}
\end{equation}%
We minimize the functional $I_{\nu ,\alpha }\left( V\right) $ on the set $S$,

\begin{equation}
V\in S=\{V\in H^{2}(\Omega ):\,V(\mathbf{x})=\psi _{0}\left( \mathbf{x}%
\right) ,\,V_{z}(\mathbf{x})=\psi _{1}\left( \mathbf{x}\right) ,\mathbf{x}%
\in \Gamma ,V(\mathbf{x})=0,\,\mathbf{x}\in \partial \Omega \setminus \Gamma
\}.  \label{eq:setW}
\end{equation}%
In (\ref{eq:Jtail}), $\alpha >0$ is the regularization parameter. The
multiplier $\exp \left( -2\mu \left( s+d\right) ^{-\nu }\right) $ is
introduced to balance two terms in the right hand side of (\ref{eq:Jtail}).

\textbf{Remark 3.2}. Since the Laplace operator is linear, one can also find
an approximate solution of problem (\ref{eq:tail}), (\ref{eq:tailbcs}) by
the regular quasi-reversibility method via setting in (\ref{eq:Jtail}) $\mu
=0$\ \cite{Klibanov2015}. However, we have noticed that a better
computational accuracy is provided in the presence of the CWF$.$\ This
observation coincides with the one of \cite{BakushinskiiKlibanov17} where it
was noticed numerically that the presence of the CWF in an analog of the
functional (\ref{eq:Jtail}) for the 1D heat equation provides a better
solution accuracy for the quasi-reversibility method.

We now follow the classical Tikhonov regularization concept \cite{T}. By
this concept, we should assume that there exists an exact solution $V_{\ast
}\left( \mathbf{x}\right) $ of the problem (\ref{eq:Jtail}), (\ref{eq:setW})
with the noiseless data $\psi _{0\ast }(\mathbf{x}),\psi _{1\ast }(\mathbf{x}%
).$ Below the subscript \textquotedblleft $\ast $" is related only to the
exact solution. In fact, however, the data $\psi _{0}(\mathbf{x}),\psi _{1}(%
\mathbf{x})$ contain noise. Let $\delta \in \left( 0,1\right) $ be the level
of noise in the data $\psi _{0}(\mathbf{x}),\psi _{1}(\mathbf{x})$. Again,
following the same concept, we should assume that the number $\delta \in
\left( 0,1\right) $ is sufficiently small. Assume that there exist functions 
$Q\left( \mathbf{x}\right) ,Q_{\ast }\left( \mathbf{x}\right) \in
H^{2}\left( \Omega \right) $ such that (see (\ref{eq:setW})) 
\begin{equation}
Q\left( \mathbf{x}\right) =\psi _{0}(\mathbf{x}),\quad \partial _{z}Q(%
\mathbf{x})=\psi _{1}(\mathbf{x}),\quad \mathbf{x}\in \Gamma ;\quad Q(%
\mathbf{x})=0,\quad \mathbf{x}\in \partial \Omega \setminus \Gamma ,
\label{3.4}
\end{equation}%
\begin{equation}
Q_{\ast }\left( \mathbf{x}\right) =\psi _{0\ast }(\mathbf{x}),\quad \partial
_{z}Q_{\ast }(\mathbf{x})=\psi _{1\ast }(\mathbf{x}),\quad \mathbf{x}\in
\Gamma ;\quad Q_{\ast }(\mathbf{x})=0,\quad \mathbf{x}\in \partial \Omega
\setminus \Gamma ,  \label{3.5}
\end{equation}%
\begin{equation}
\left\Vert Q-Q_{\ast }\right\Vert _{H^{3}\left( \Omega \right) }<\delta .
\label{3.6}
\end{equation}%
Introduce the number $t_{\nu },$%
\begin{equation}
t_{\nu }=\left( s-\xi \right) ^{-\nu }-\left( s+d\right) ^{-\nu }>0.
\label{3.70}
\end{equation}

Let 
\begin{equation}
W\left( \mathbf{x}\right) =V\left( \mathbf{x}\right) -Q\left( \mathbf{x}%
\right) .  \label{1}
\end{equation}
Then by (\ref{eq:Jtail}) and (\ref{eq:setW}) the functional $I_{\mu ,\alpha }
$ becomes 
\begin{equation}
\widetilde{I}_{\mu ,\alpha }\left( W\right) =\exp \left( -2\mu \left(
s+d\right) ^{-\nu }\right) \int_{\Omega }\left\vert \Delta W+\Delta
Q\right\vert ^{2}\varphi _{\mu }\left( z\right) d\mathbf{x}+\alpha \Vert
W+Q\Vert _{H^{3}(\Omega )}^{2},\text{ }W\in H_{0}^{3}\left( \Omega \right) .
\label{5.1}
\end{equation}%
Theorem \ref{thm:4.2} of section 4 claims that for each $\alpha >0$ there exists
unique minimizer $W_{\mu ,\nu ,\alpha }\in H^{3}\left( \Omega \right) $ of
the functional (\ref{eq:Jtail}), which is called the \textquotedblleft
regularized solution". Using (\ref{1}), denote $V_{\mu ,\nu ,\alpha }=W_{\mu
,\nu ,\alpha }+Q.$ It is stated in Theorem \ref{thm:4.2} that one can choose a
sufficiently large number $\nu _{0}=\nu _{0}\left( \Omega ,s\right) $
depending only on $\Omega $ and $s$ such that for any fixed value of the
parameter $\nu \geq \nu _{0}$ the choices 
\begin{equation}
\alpha =\alpha \left( \delta \right) =\delta ,\mu =\ln \left( \delta
^{-1/\left( 2t_{\nu }\right) }\right)   \label{3.71}
\end{equation}%
regularized solutions converge to the exact solution as $\delta \rightarrow
0.$ More precisely, there exists a constant $C=C\left( \Omega \right) >0$%
\emph{\ }such that 
\begin{equation}
\left\Vert V_{\mu \left( \delta \right) ,\nu ,\alpha \left( \delta \right)
}-V_{\ast }\right\Vert _{H^{2}\left( \Omega \right) }\leq C\left( 1+\Vert
V_{\ast }\Vert _{H^{3}(\Omega )}\right) \sqrt{\delta }\sqrt{\ln \left(
\delta ^{-1/\left( 2t_{\nu }\right) }\right) }.  \label{3.7}
\end{equation}%
Here and below $C=C\left( \Omega \right) >0$ denotes different positive
constants depending only on the domain $\Omega .$

\subsection{Associated spaces}

\label{sec:3.4}

Below, for any complex number $z\in \mathbb{C}$ we denote $\overline{z}$ its
complex conjugate. It is convenient for us to consider any complex valued
function $U=\mathop{\rm Re}U+i\mathop{\rm
Im}U=U_{1}+iU_{2}$ as the 2D vector function $U=\left( U_{1},U_{2}\right) .$
Thus, below any Banach space we use for a complex valued function is
actually the space of such 2D real valued vector functions. Norms in these
spaces of 2D vector functions are defined in the standard way, so as scalar
products, in the case of Hilbert spaces. For brevity we do not differentiate
below between complex valued functions and corresponding 2D vector
functions. However, it is always clear from the context what is what.

We define the Hilbert space $H_{m}$ of complex valued functions $f\left( \mathbf{x},k\right) $ as%
\begin{equation}
H_{m}=\left\{ f\left( \mathbf{x},k\right) :\left\Vert f\right\Vert _{H_{m}}=%
\left[ \int_{\underline{k}}^{\overline{k}}\left\Vert f\left( \mathbf{x%
},k\right) \right\Vert _{H^{m}\left( \Omega \right) }^{2}dk\right]
^{1/2}<\infty \right\} ,\text{ }m=1,2,3.  \label{3.700}
\end{equation}%
Denote $\left[ ,\right] $ the scalar product in the space $H_{3}.$ The
subspace $H_{m}^{0}$ of the space $H_{m}$ is defined as%
\begin{equation*}
H_{m}^{0}=\left\{ f\in H_{m}:f\left( \mathbf{x},k\right) \mid _{\partial
\Omega }=0,f_{z}\left( \mathbf{x},k\right) \mid _{\Gamma }=0,\forall k\in %
\left[ \underline{k},\overline{k}\right] \right\} .
\end{equation*}%
Also, in the case of functions independent on $k$,%
\begin{equation*}
H_{0}^{m}\left( \Omega \right) =\left\{ f\left( \mathbf{x}\right) \in
H^{m}\left( \Omega \right) :f\left( \mathbf{x}\right) \mid _{\partial \Omega
}=0,f_{z}\left( \mathbf{x}\right) \mid _{\Gamma }=0\right\} .
\end{equation*}%
Similarly we define for $r=0,1,2$%
\begin{equation*}
C_{r}=\left\{ f\left( \mathbf{x},k\right) :\left\Vert f\right\Vert
_{C_{r}}=\max_{k\in \left[ \underline{k},\overline{k}\right] }\left\Vert
f\left( \mathbf{x},k\right) \right\Vert _{C^{r}\left( \overline{\Omega }%
\right) }\right\} ,
\end{equation*}%
where $C^{0}\left( \overline{\Omega }\right) =C\left( \overline{\Omega }%
\right) .$ Embedding theorem implies that: 
\begin{equation}
H_{3+r}\subset C_{1+r},\left\Vert f\right\Vert _{C_{1+r}}\leq C\left\Vert
f\right\Vert _{H_{3+r}},\text{ }\forall f\in H_{3+r},r=0,1\text{,}
\label{3.73}
\end{equation}%
\begin{equation}
\left\Vert \widetilde{f}\right\Vert _{C^{1}\left( \overline{\Omega }\right)
}\leq C\left\Vert \widetilde{f}\right\Vert _{H^{3}\left( \Omega \right) },%
\text{ }\forall \widetilde{f}\in H^{3}\left( \Omega \right) .  \label{3.74}
\end{equation}

\subsection{The weighted Tikhonov-like functional}

\label{sec:3.5}

Suppose that there exists a function $F\left( \mathbf{x},k\right) \in H_{4}$
such that (see (\ref{eq:intdiffbcs})): 
\begin{equation}
F\left( \mathbf{x},k\right) \mid _{\Gamma }=\phi _{0}\left( \mathbf{x}%
,k\right) ,\text{ }F_{z}\left( \mathbf{x},k\right) \mid _{\Gamma }=\phi
_{1}\left( \mathbf{x},k\right) ,\text{ }F\left( \mathbf{x},k\right) \mid
_{\partial \Omega \diagdown \Gamma }=0.  \label{3.8}
\end{equation}%
Also, assume that there exists an exact solution $c_{\ast }\left( \mathbf{x}\right) $ of our CISP satisfying the above conditions imposed on the
coefficient $c\left( \mathbf{x}\right) $ and generating the noiseless
boundary data $\phi _{0\ast }$ and $\phi _{1\ast }$ in (\ref{eq:intdiffbcs}%
). Let the function $F_{\ast }\left( \mathbf{x},k\right) \in H_{3}$
satisfies boundary conditions (\ref{3.8}) in which functions $\phi _{0}$ and 
$\phi _{1}$ are replaced with functions $\phi _{0\ast }$ and $\phi _{1\ast }$
respectively. We assume that%
\begin{equation}
\left\Vert F-F_{\ast }\right\Vert _{H_{4}}<\delta .  \label{3.80}
\end{equation}

Let $q_{\ast }\in H_{3}$ be the function $q$ generated by the exact
coefficient $c_{\ast }\left( \mathbf{x}\right) .$ Introduce functions $%
p,p_{\ast }\in H_{3}^{0}$ as 
\begin{equation}
p\left( \mathbf{x},k\right) =q\left( \mathbf{x},k\right) -F\left( \mathbf{x}%
,k\right) ,\text{ }p_{\ast }\left( \mathbf{x},k\right) =q_{\ast }\left( 
\mathbf{x},k\right) -F_{\ast }\left( \mathbf{x},k\right) .  \label{3.9}
\end{equation}%
It follows from the discussion in section 2.2 about the smoothness as well
as from (\ref{eq:v}), (\ref{eq:q}) and (\ref{3.9}) that the functions $%
p,p_{\ast }\in H_{3}^{0}.$ Let $R>0$ be an arbitrary number. Consider the
ball $B\left( R\right) \subset H_{3}^{0}$ of the radius $R$,%
\begin{equation}
B\left( R\right) =\left\{ f\in H_{3}^{0}:\left\Vert f\right\Vert
_{H_{3}}<R\right\} .  \label{3.10}
\end{equation}

Based on the integro-differential equation (\ref{eq:intdiff}), boundary
conditions (\ref{eq:intdiffbcs}) for it, (\ref{3.8}) and (\ref{3.9}), we
construct our weighted Tikhonov-like functional with the CWF (\ref{3.3}) in
it as 
\begin{equation}
J_{\lambda ,\rho }\left( p\right) =\exp \left( -2\lambda \left( s+d\right)
^{-\nu }\right) \int_{\underline{k}}^{\overline{k}%
}\int_{\Omega }\left\vert L\left( p+F\right) \left( \mathbf{x},\kappa
\right) \right\vert ^{2}\varphi _{\lambda }^{2}\left( z\right) d\mathbf{x}%
d\kappa +\rho \left\Vert p\right\Vert _{H_{3}}^{2},  \label{eq:J}
\end{equation}%
where $\rho >0$ is the regularization parameter. Similarly with (\ref%
{eq:Jtail}), the multiplier $\exp \left( -2\lambda \left( s+d\right) ^{-\nu
}\right) $ is introduced to balance two terms in the right hand side of (\ref%
{eq:J}). The minimizer $V_{\mu \left( \delta \right) ,\nu ,\alpha \left(
\delta \right) }$ of the functional (\ref{eq:Jtail}) is chosen in $%
J_{\lambda ,\rho }\left( p\right) $ as the tail function. We consider the
following minimization problem:

\bigskip \textbf{Minimization Problem}. \emph{Minimize the functional }$J$%
\emph{$_{\lambda ,\rho }$}$(q)$\emph{\ on the set $\overline{B\left(
R\right) }$.}

\section{Theorems}

\label{sec:4}

In this section we formulate theorems about numerical procedures considered
in section 3. We start from the Carleman estimate with the CWF (\ref{3.3}).

 \begin{theorem}[Carleman estimate]
 \label{thm:4.1}
 Let $\Omega \subset  \mathbb{R}^{3}$ be the above domain (\ref{eq:2.1}). Temporary denote $\mathbf{x}=\left( x,y,z\right) =\left( x_{1},x_{2},x_{3}\right)$. There exist numbers $C=C\left( \Omega \right) >0$, $\nu _{0}=\nu
_{0}\left( \Omega ,s,d\right) \geq 1$ and $\lambda _{0}=\lambda
_{0}\left( \Omega ,s,d\right) \geq 1$ depending only on listed
parameters such that for any real valued function $u\in H_{0}^{2}\left(
\Omega \right) $ the following Carleman estimate holds with the CWF $\varphi _{\lambda }\left( z\right) $ in (\ref{3.3} for and
fixed number $\nu \geq \nu _{0}$ and for all $\lambda \geq \lambda _{0}$
\begin{equation}
\int_{\Omega }\left( \Delta u\right) ^{2}\varphi _{\lambda }\left(
z\right) d\mathbf{x}\geq \frac{C}{\lambda }\sum_{i,j=1}^{3}\int_{\Omega }\left( u_{x_{i}x_{j}}\right) ^{2}\varphi _{\lambda }\left(
z\right) d\mathbf{x}+C\lambda \int_{\Omega }\left( \nabla u\right)
^{2}\varphi _{\lambda }\left( z\right) d\mathbf{x}+C\lambda
^{3}\int_{\Omega }u^{2}\varphi _{\lambda }\left( z\right) d\mathbf{x}
\label{4.1}
\end{equation}
\end{theorem}

\textbf{Remark 4.1}. A close analog of Theorem \ref{thm:4.1} is formulated as lemma
4.1 of \cite{KlibanovThanh15} and is proven in the proof of lemma 6.5.1 of 
\cite{BeilinaKlibanov12}. Hence, we omit the proof of Theorem \ref{thm:4.1}.

The next theorem is about the problem (\ref{eq:Jtail}), (\ref{eq:setW}).

 \begin{theorem}
 \label{thm:4.2}
Assume that there exists a function $Q\in H^{3}\left( \Omega \right) $ satisfying conditions (\ref{3.4}), (\ref{3.6}). Then for each set of parameters $\mu ,\nu ,\alpha >0$ there
exists unique minimizer $W_{\mu ,\alpha ,\nu }\in H^{3}\left( \Omega\right) $ of the functional (\ref{5.1}). Let $V_{\mu ,\nu ,\alpha}=W_{\mu ,\nu ,\alpha }+Q$ (see (\ref{1})). Suppose that there
exists an exact solution $V_{\ast }\in H^{3}\left( \Omega \right) $
of the problem (\ref{eq:tail}), (\ref{eq:tailbcs}) with the noiseless
boundary data $\psi _{0\ast }(x),\psi _{1\ast }(x)$. Also, assume
that there exists a function $Q_{\ast }\in H^{3}\left( \Omega \right) $%
  satisfying conditions (\ref{3.5}) and such that 
\begin{equation}
\left\Vert Q_{\ast }\right\Vert _{H^{3}\left( \Omega \right) }\emph{\ }\leq
C\left\Vert V_{\ast }\right\Vert _{H^{3}\left( \Omega \right) }.
\label{4.200}
\end{equation}%
Let inequality (\ref{3.6}) hold, where $\delta \in \left( 0,1\right) $ is the level of noise in the data. Let $\nu _{0}\left( \Omega,s\right) $ and $\lambda _{0}\left( \Omega ,s\right) $ becnumbers of Theorem \ref{thm:4.1}. Fix a number $\nu \geq \nu _{0}\left( \Omega,s\right) $ and let the parameter $\nu $  be independent on $\delta$. Choose a number $\delta _{0}\in \left( 0,e^{-2\lambda_{0}t_{\nu }}\right)$, where $\lambda _{0}$  is defined in
Theorem \ref{thm:4.1} and the number $t_{\nu }$  is defined in (\ref{3.70}).
For any $\delta \in \left( 0,\delta _{0}\right) $  let the choice (\ref{3.71}) holds. Then the convergence estimate (\ref{3.7}) of functions $V_{\mu \left( \delta \right) ,\nu ,\alpha \left( \delta \right) }$  to the exact solution $V_{\ast }$  holds for $\delta \rightarrow 0$. In addition, the function $V_{\mu \left( \delta \right) ,\nu ,\alpha \left( \delta \right) }\in C^{1}\left( \overline{\Omega }\right) $ and
\begin{equation}
C\left\Vert V_{\mu \left( \delta \right) ,\nu ,\alpha \left( \delta \right)
}\right\Vert _{C^{1}\left( \overline{\Omega }\right) }\leq \left\Vert V_{\mu
\left( \delta \right) ,\nu ,\alpha \left( \delta \right) }\right\Vert
_{H^{3}\left( \Omega \right) }\leq C\left( 1+\left\Vert V_{\ast }\right\Vert
_{H^{3}\left( \Omega \right) }\right) .  \label{4.100}
\end{equation}
 \end{theorem}
Theorem \ref{thm:4.3} is the central analytical result of this paper.

\begin{theorem}[Global strict convexity]
\label{thm:4.3}
Assume that conditions of Theorem \ref{thm:4.2} hold. Set in (\ref{eq:intdiff}) $V=V_{\mu \left(\delta \right) ,\nu ,\alpha \left( \delta \right) },$  where the
function $V_{\mu \left( \delta \right) ,\nu ,\alpha \left( \delta \right) }$%
 is the one of Theorem \ref{thm:4.2}. The functional $J_{\lambda ,\rho }\left(
p\right) $  has the Fr\'{e}chet derivative $J_{\lambda ,\rho
}^{\prime }\left( p\right) \in H_{3}$ at any point $p\in H_{3}^{0}.$%
  Assume that there exists a function $F,F_{\ast }\left( \mathbf{x}%
,k\right) \in H_{4}$  satisfying conditions (\ref{3.8}), (\ref{3.80}),
where $\delta \in \left( 0,1\right) .$ Let $\lambda _{0}=\lambda
_{0}\left( \Omega \right) $  be the number defined in Theorem \ref{thm:4.1}. Then
there exists a number $\lambda _{1}=\lambda _{1}\left( \Omega ,R,\left\Vert
F_{\ast }\right\Vert _{H_{4}},\left\Vert V_{\ast }\right\Vert _{H^{3}\left(
\Omega \right) },\underline{k},\overline{k}\right) \geq \lambda _{0}\left(
\Omega \right) $  and a number $C_{1}=C_{1}\left( \Omega
,R,\left\Vert F_{\ast }\right\Vert _{H_{4}},\left\Vert V_{\ast }\right\Vert
_{H^{3}\left( \Omega \right) },\underline{k},\overline{k}\right) >0,$ 
both depending only on listed parameters, such that for any $\lambda
\geq \lambda _{1}$  the functional $J_{\lambda ,\rho }\left( p\right) 
$  is strictly convex on $\overline{B\left( R\right) }.$ In
other words, the following estimates are valid for all $p_{1},p_{2}\in 
\overline{B\left( R\right) }:$
\begin{equation}
J_{\lambda ,\rho }\left( p_{2}\right) -J_{\lambda ,\rho }\left( p_{1}\right)
-J_{\lambda ,\rho }^{\prime }\left( p_{1}\right) \left( p_{2}-p_{1}\right)
\geq \frac{C_{1}}{\lambda }\left\Vert p_{2}-p_{1}\right\Vert
_{H_{2}}^{2}+\rho \left\Vert p_{2}-p_{1}\right\Vert _{H_{3}}^{2},
\label{4.2}
\end{equation}%
\begin{equation}
J_{\lambda ,\rho }\left( p_{2}\right) -J_{\lambda ,\rho }\left( p_{1}\right)
-J_{\lambda ,\rho }^{\prime }\left( p_{1}\right) \left( p_{2}-p_{1}\right)
\geq C_{1}\left\Vert p_{2}-p_{1}\right\Vert _{H_{1}}^{2}+\rho \left\Vert
p_{2}-p_{1}\right\Vert _{H_{3}}^{2}.  \label{4.3}
\end{equation}
\end{theorem}
\textbf{Remark 4.2}. The first term in the right hand side of (\ref{4.3})
does not decay with the increase of $\lambda ,$ unlike (\ref{4.2}). Hence,
the \textquotedblleft convexity property" of the functional $J_{\lambda
,\rho }$ is sort of better in terms of the $H_{1}-$norm in (\ref{4.3})
rather than in terms of the $H_{2}-$norm in (\ref{4.2}). On the other hand,
the norm of that term is weaker than the one in (\ref{4.2}). Also, to
establish convergence of reconstructed coefficients $c_{n}\left( \mathbf{x}%
\right) ,$ we need the $H_{2}-$norm: see (\ref{4.11}) and Remark 3.1.

\begin{theorem}
\label{thm:4.4}
Suppose that the conditions of Theorems \ref{thm:4.2} and
\ref{thm:4.3} regarding the tail function $V=V_{\mu \left( \delta \right) ,\nu
,\alpha \left( \delta \right) }$  and the function $F$ hold.
The Fr\'{e}chet derivative $J_{\lambda ,\rho }^{\prime }$ of the
functional $J_{\lambda ,\rho }$  satisfies the Lipschitz continuity
condition in any ball $B\left( R^{\prime }\right) $ as in (\ref{3.10}%
) with any $R^{\prime }>0.$ In other words, the following inequality
holds with the constant $M=M\left( \Omega ,R^{\prime },\left\Vert F_{\ast
}\right\Vert _{H_{4}},\left\Vert V_{\ast }\right\Vert _{H^{3}\left( \Omega
\right) },\lambda ,\nu ,\rho ,\underline{k},\overline{k}\right) >0$
depending only on listed parameters:
\begin{equation*}
\left\Vert J_{\lambda ,\rho }^{\prime }\left( p_{1}\right) -J_{\lambda ,\rho
}^{\prime }\left( p_{2}\right) \right\Vert _{H_{3}}\leq M\left\Vert
p_{1}-p_{2}\right\Vert _{H_{3}},\text{ }\forall p_{1},p_{2}\in B\left(
R^{\prime }\right) .
\end{equation*}
\end{theorem}

Let $P_{\overline{B}}:H_{3}^{0}\rightarrow \overline{B\left( R\right) }$ be
the projection operator of the Hilbert space $H_{3}^{0}$ on $\overline{%
B\left( R\right) }.$ Let $p_{0}\in B\left( R\right) $ be an arbitrary point
of the ball $B\left( R\right) $. Consider the following sequence: 
\begin{equation}
p_{n}=P_{\overline{B}}\left( p_{n-1}-\omega J_{\lambda ,\rho }^{\prime
}\left( p_{n-1}\right) \right) ,\text{ }n=1,2,...,  \label{4.9}
\end{equation}%
where $\omega \in \left( 0,1\right) $ is a certain number.

\begin{theorem}
\label{thm:4.5}
Assume that conditions of Theorems \ref{thm:4.2} and \ref{thm:4.3}
hold. Let $\lambda \geq \lambda _{1},$ where $\lambda _{1}$ is the number of Theorem \ref{thm:4.3}. Then there exists unique minimizer $p_{\min
,\lambda }\in \overline{B\left( R\right) }$  of the functional $%
J_{\lambda ,\rho }\left( p\right) $  on the set $\overline{B\left(
R\right) }$ and
\begin{equation}
J_{\lambda ,\rho }^{\prime }\left( p_{\min ,\lambda }\right) \left(
y-p_{\min ,\lambda }\right) \geq 0,\text{ \ }\forall y\in H_{3}^{0}.
\label{4.6}
\end{equation}%
Also, there exists a sufficiently small number $\omega _{0}=\omega
_{0}\left( \Omega ,R,\left\Vert F\right\Vert _{H_{4}},\left\Vert V_{\ast
}\right\Vert _{H^{3}\left( \Omega \right) },\underline{k},\overline{k}%
,\lambda ,\delta \right) \in \left( 0,1\right) $  depending only on
listed parameters such that for any $\omega \in \left( 0,\omega _{0}\right) 
$  the sequence (\ref{4.9}) converges to the minimizer $p_{\min
,\lambda }\in \overline{B\left( R\right) }$  of the functional $%
J_{\lambda ,\rho }\left( p\right) $  on the set $\overline{B\left(
R\right) }$,
\begin{equation}
\left\Vert p_{\min ,\lambda }-p_{n}\right\Vert _{H_{3}}\leq r^{n}\left\Vert
p_{\min ,\lambda }-p_{0}\right\Vert _{H_{3}},\text{ }n=1,2,..,  \label{4.90}
\end{equation}%
 where the number $r=r\left( \omega ,\Omega ,R,\left\Vert
F\right\Vert _{H_{4}},\left\Vert V_{\ast }\right\Vert _{H^{3}\left( \Omega
\right) },\underline{k},\overline{k},\lambda ,\delta \right) \in \left(
0,1\right) $ depends only on listed parameters.
\end{theorem}

By (\ref{4.90}) we estimate the convergence rate of the sequence (\ref{4.9})
to the minimizer. The next question is about the convergence of this
sequence to the exact solution $p_{\ast }$ assuming that it exists.

\begin{theorem}
\label{thm:4.6}
Assume that conditions of Theorems \ref{thm:4.2} and \ref{thm:4.3}
hold. Let $\lambda _{1}$  be the number of Theorem \ref{thm:4.3}. Choose a
number $\delta _{1}\in \left( 0,e^{-2\lambda _{1}t_{\nu }}\right) .$ For $%
\delta \in \left( 0,\delta _{1}\right) ,$ set $\rho =\rho \left(
\delta \right) =\sqrt{\delta },\lambda =\lambda \left( \delta \right) =\ln
\left( \delta ^{-1/\left( 2t_{\nu }\right) }\right) .$ Furthermore,
assume that the exact solution $p_{\ast }$ exists and $p_{\ast }\in
B\left( R\right) $. Then there exists a number $C_{2}=C_{2}\left(
\Omega ,R,\left\Vert F\right\Vert _{H_{4}},\left\Vert V_{\ast }\right\Vert
_{H^{3}\left( \Omega \right) },\underline{k},\overline{k}\right) >0$
depending only on listed parameters such that
\begin{equation}
\left\Vert p_{\ast }-p_{\min ,\lambda \left( \delta \right) }\right\Vert
_{H_{2}}\leq C_{2}\delta ^{1/4}\left[ \ln \left( \delta ^{-1/\left( 2t_{\nu
}\right) }\right) \right] ^{3/4},  \label{4.7}
\end{equation}%
\begin{equation}
\left\Vert c_{\ast }-c_{\min ,\lambda \left( \delta \right) }\right\Vert
_{L_{2}\left( \Omega \right) }\leq C_{2}\delta ^{1/4}\left[ \ln \left(
\delta ^{-1/\left( 2t_{\nu }\right) }\right) \right] ^{3/4},  \label{4.8}
\end{equation}%
where the function $c_{\min ,\lambda \left( \delta \right) }\left( 
\mathbf{x}\right) $ is reconstructed from the function $p_{\min ,\lambda
\left( \delta \right) }$ using (\ref{3.9}) and Remark 3.1. In addition, the following convergence estimates hold
\begin{equation}
\left\Vert p_{\ast }-p_{n}\right\Vert _{H_{2}}\leq C_{2}\delta ^{1/4}\left[
\ln \left( \delta ^{-1/\left( 2t_{\nu }\right) }\right) \right]
^{3/4}+r^{n}\left\Vert p_{\min ,\lambda \left( \delta \right)
}-p_{0}\right\Vert _{H_{3}},\text{ }n=1,2,...,  \label{4.10}
\end{equation}%
\begin{equation}
\left\Vert c_{\ast }-c_{n}\right\Vert _{L_{2}\left( \Omega \right) }\leq
C_{2}\delta ^{1/4}\left[ \ln \left( \delta ^{-1/\left( 2t_{\nu }\right)
}\right) \right] ^{3/4}+C_{2}r^{n}\left\Vert p_{\min ,\lambda \left( \delta
\right) }-p_{0}\right\Vert _{H_{3}},\text{ }n=1,2,...,  \label{4.11}
\end{equation}%
where $r$ is the number in (\ref{4.90}) and the function $c_{n}\left( \mathbf{x}\right) $  is reconstructed from the function $p_{n}\left( \mathbf{x},k\right) $ using (\ref{3.9}) and Remark 3.1.
\end{theorem}

\textbf{Remark 4.1}. Since $R>0$\ is an arbitrary number and $p_{0}$\ is an
arbitrary point of the ball $B\left( R\right) $, then Theorems \ref{thm:4.5} and \ref{thm:4.6}
ensure the global convergence of the gradient projection method for our
case, see the second paragraph of section 1. We note that if a functional is
non convex, then the convergence of a gradient-like method of its
minimization might be guaranteed only if the starting point of iterations is
located in a sufficiently small neighborhood of its minimizer.

\section{Proofs}

\label{sec:5}

In this section we prove theorems formulated in section 4, except of Theorem
\ref{thm:4.1} (see Remark 4.1).

\subsection{Proof of Theorem \ref{thm:4.2}}

\label{sec:5.1}

By (\ref{1}) and (\ref{5.1}) the vector function $W_{\min }=\left( W_{1,\min
},W_{2,\min }\right) \in H_{0}^{3}\left( \Omega \right) $ is a minimizer of
the functional $\widetilde{I}_{\mu ,\alpha }\left( W\right) $ if and only if 
\begin{equation}
\begin{gathered} \exp \left( -2\mu \left( s+d\right) ^{-\nu }\right)
\int_{\Omega }\left( \Delta W_{1,\min }\Delta h_{1}+\Delta W_{2,\min
}\Delta h_{2}\right) \varphi _{\mu }\left( z\right) d\mathbf{x}+\alpha
\left( \left( W_{\min },h\right) \right) = \\ - \exp \left( -2\mu \left(
s+d\right) ^{-\nu }\right) \int_{\Omega }\left( \Delta Q_{1}\Delta
h_{1}+\Delta Q_{2}\Delta h_{2}\right) \varphi _{\mu }\left( z\right)
d\mathbf{x}-\alpha \left( \left( Q,h\right) \right) ,\,\forall h=\left(
h_{1},h_{2}\right) \in H_{0}^{3}\left( \Omega \right) , \end{gathered}
\label{5.2}
\end{equation}%
where $\left( \left( ,\right) \right) $ is the scalar product in $%
H^{3}\left( \Omega \right) .$ For any vector function $P=\left(
P_{1},P_{2}\right) \in H_{0}^{3}\left( \Omega \right) $ consider the
expression in the left hand side of (\ref{5.2}) in which the vector $\left(
W_{1,\min },W_{2,\min }\right) $ is replaced with $\left( P_{1},P_{2}\right)
.$ Then this expression defines a new scalar product $\left\{ P,h\right\} $
in $H^{3}\left( \Omega \right) ,$\ and the corresponding norm $\sqrt{\left\{
P,P\right\} }$ is equivalent to the norm in $H^{3}\left( \Omega \right) .$
Next, 
\begin{equation*}
\begin{gathered} \left\vert -\exp \left( -2\mu \left( s+d\right) ^{-\nu
}\right) \int_{\Omega }\left( \Delta Q_{1}\Delta h_{1}+\Delta
Q_{2}\Delta h_{2}\right) \varphi _{\mu }\left( z\right) d\mathbf{x}-\alpha
\left( \left( Q,h\right) \right) \right\vert \leq D\left\Vert Q\right\Vert
_{H^{3}\left( \Omega \right) }\left\Vert h\right\Vert _{H^{3}\left( \Omega
\right) }\\ \leq D_{1}\sqrt{\left\{ Q,Q\right\} }\sqrt{\left\{ h,h\right\}
},\text{ }\forall h=\left( h_{1},h_{2}\right) \in H_{0}^{3}\left( \Omega
\right) \end{gathered}
\end{equation*}%
with certain constants $D,D_{1}$ independent on $Q$ and $h$ but dependent on
parameters $\mu ,\nu .$ Hence, Riesz theorem implies that there exists
unique vector function $\widehat{Q}=\widehat{Q}\left( Q\right) \in
H_{0}^{3}\left( \Omega \right) $ such that 
\begin{equation*}
-\exp \left( -2\mu \left( s+d\right) ^{-\nu }\right) \int_{\Omega
}\left( \Delta Q_{1}\Delta h_{1}+\Delta Q_{2}\Delta h_{2}\right) \varphi
_{\mu }\left( z\right) d\mathbf{x}-\alpha \left( \left( Q,h\right) \right)
=\left\{ \widehat{Q},h\right\} ,\text{ }\forall h=\left( h_{1},h_{2}\right)
\in H_{0}^{3}\left( \Omega \right) .
\end{equation*}%
Hence, by (\ref{5.2}) $\left\{ W_{\min },h\right\} =\left\{ \widehat{Q}%
,h\right\} ,\forall h\in H_{0}^{3}\left( \Omega \right) .$ Hence, $W_{\min }=%
\widehat{Q}.$ Thus, existence and uniqueness of the minimizer of the
functional $\widetilde{I}_{\mu ,\alpha }\left( W\right) $ are established,
and the same for $I_{\mu ,\alpha }\left( V\right) $.

We now prove convergence estimate (\ref{3.7}). Let $W_{\ast }=V_{\ast
}-Q_{\ast }\in H_{0}^{3}\left( \Omega \right) .$ Denote $\widetilde{W}%
=W_{\min }-W_{\ast },$ $\widetilde{Q}=Q-Q_{\ast }.$ Since%
\begin{equation}
\exp \left( -2\mu \left( s+d\right) ^{-\nu }\right) \int_{\Omega
}\left( \Delta W_{\ast ,1}\Delta h_{1}+\Delta W_{\ast ,2}\Delta h_{2}\right)
\varphi _{\mu }\left( z\right) d\mathbf{x}+\alpha \left[ W_{\ast },h\right] 
\label{5.3}
\end{equation}%
\begin{equation*}
=-\exp \left( -2\mu \left( s+d\right) ^{-\nu }\right) \int_{\Omega
}\left( \Delta Q_{1}^{\ast }\Delta h_{1}+\Delta Q_{2}^{\ast }\Delta
h_{2}\right) \varphi _{\mu }\left( z\right) d\mathbf{x+}\alpha \left[
W_{\ast },h\right] ,\text{ }\forall h\in H_{0}^{3}\left( \Omega \right) ,
\end{equation*}%
then subtracting (\ref{5.3}) from (\ref{5.2}) and setting $h=\widetilde{W},$
we obtain%
\begin{equation*}
\begin{gathered} \exp \left( -2\mu \left( s+d\right) ^{-\nu }\right)
\int_{\Omega }\left( \Delta \widetilde{W}\right) ^{2}\varphi _{\mu
}\left( z\right) d\mathbf{x+}\alpha \left\Vert \widetilde{W}\right\Vert
_{H^{3}\left( \Omega \right) }^{2}\\ =-\exp \left( -2\mu \left( s+d\right)
^{-\nu }\right) \int_{\Omega }\left( \Delta \widetilde{Q}_{1}\Delta
\widetilde{W}_{1}+\Delta \widetilde{Q}_{2}\Delta \widetilde{W}_{2}\right)
\varphi _{\mu }\left( z\right) d\mathbf{x-}\alpha \left( \left( W_{\ast
}+Q,\widetilde{W}\right) \right) . \end{gathered}
\end{equation*}%
Using the Cauchy-Schwarz inequality, taking into account (\ref{3.6}) and
recalling that $\alpha =\delta $, we obtain%
\begin{equation}
\begin{gathered} \exp \left( -2\mu \left( s+d\right) ^{-\nu }\right)
\int_{\Omega }\left( \Delta \widetilde{W}\right) ^{2}\varphi _{\mu
}\left( z\right) d\mathbf{x}+\delta \left\Vert \widetilde{W}\right\Vert
_{H^{3}\left( \Omega \right) }^{2} \\ \leq C\delta \left( 1+\left\Vert
V_{\ast }\right\Vert _{H^{3}\left( \Omega \right) }^{2}\right) +C\exp \left(
2\mu t_{\nu }\right) \delta ^{2}. \end{gathered}  \label{5.4}
\end{equation}%
Since $\mu =\ln \left( \delta ^{-1/\left( 2t_{\nu }\right) }\right) $ and $%
\delta \in \left( 0,1\right) ,$ then $\exp \left( 2\mu t_{\nu }\right)
\delta ^{2}=\delta $ and 
\begin{equation*}
C\delta \left( 1+\left\Vert V_{\ast }\right\Vert _{H^{3}\left( \Omega
\right) }^{2}\right) +C\exp \left( 2\mu t_{\nu }\right) \delta ^{2}\leq
C\delta \left( 1+\left\Vert V_{\ast }\right\Vert _{H^{3}\left( \Omega
\right) }^{2}\right) ,
\end{equation*}%
then (\ref{5.4}) implies that%
\begin{equation}
\left\Vert \widetilde{W}\right\Vert _{H^{3}\left( \Omega \right) }\leq
C\left( 1+\left\Vert V_{\ast }\right\Vert _{H^{3}\left( \Omega \right)
}\right) ,  \label{5.5}
\end{equation}%
\begin{equation}
\exp \left( -2\mu \left( s+d\right) ^{-\nu }\right) \int_{\Omega
}\left( \Delta \widetilde{W}\right) ^{2}\varphi _{\mu }\left( z\right) d%
\mathbf{x\leq }C\delta \left( 1+\left\Vert V_{\ast }\right\Vert
_{H^{3}\left( \Omega \right) }^{2}\right) .  \label{5.50}
\end{equation}%
Since 
\begin{equation*}
\exp \left( -2\mu \left( s+d\right) ^{-\nu }\right) \varphi _{\mu }\left(
z\right) \geq \exp \left( -2\mu \left( s+d\right) ^{-\nu }\right) \exp
\left( 2\mu \left( s+d\right) ^{-\nu }\right) =1,
\end{equation*}%
then Theorem \ref{thm:4.1} implies that%
\begin{equation}
\begin{gathered} \exp \left( -2\mu \left( s+d\right) ^{-\nu }\right)
\int_{\Omega }\left( \Delta \widetilde{W}\right) ^{2}\varphi _{\mu
}\left( z\right) d\mathbf{x}\\ \mathbf{\geq }\frac{C}{\mu }\left(
\sum_{i,j=1}^{3}\int_{\Omega
}\widetilde{W}_{x_{i}x_{j}}^{2}d\mathbf{x+}\mu ^{2}\int_{\Omega
}\left( \left( \nabla \widetilde{W}\right) ^{2}+\widetilde{W}^{2}\right)
d\mathbf{x}\right) \geq \frac{C}{\mu }\left\Vert \widetilde{W}\right\Vert
_{H^{2}\left( \Omega \right) }^{2}. \end{gathered}  \label{5.51}
\end{equation}%
The right estimate (\ref{4.100}) follows from (\ref{5.5}), (\ref{3.6}) and (%
\ref{4.200}). The left estimate (\ref{4.100}) follows from (\ref{3.74}).
Comparing (\ref{5.50}) with (\ref{5.51}) and recalling (\ref{1}) and (\ref%
{4.200}), we obtain (\ref{3.7}). $\square $

\subsection{Proof of Theorem \ref{thm:4.3}}

\label{sec:5.2}

Recall that we treat any complex valued function $U=\mathop{\rm Re}U+i%
\mathop{\rm Im}U=U_{1}+iU_{2}$ in two ways: (1) in its original complex
valued form and (2) in an equivalent form as a 2D vector function $\left(
U_{1},U_{2}\right) $ (section 3.4). It is always clear from the content what
is what.

Let two arbitrary functions $p_{1},p_{2}\in \overline{B\left( R\right) }.$
Denote $h=p_{2}-p_{1}.$ Then $h=\left( h_{1},h_{2}\right) \in
H_{0}^{3}\left( \Omega \right) .$ In this proof $C_{1}=C_{1}\left( \Omega
,R,\left\Vert F_{\ast }\right\Vert _{H_{4}},\left\Vert V_{\ast }\right\Vert
_{H^{3}\left( \Omega \right) },\underline{k},\overline{k}\right) >0$ denotes
different positive constants. Also, in this proof we denote for brevity $%
V\left( \mathbf{x}\right) =$\emph{\ }$V_{\mu \left( \delta \right) ,\nu
,\alpha \left( \delta \right) }\left( \mathbf{x}\right) .$ We note that due
to (\ref{3.73}), (\ref{3.74}), (\ref{3.80}), (\ref{3.10}) and (\ref{4.100})%
\begin{equation}
\left\Vert \nabla V\right\Vert _{C\left( \overline{\Omega }\right)
},\left\Vert F\right\Vert _{C^{2}\left( \overline{\Omega }\right) }\leq
C_{1},  \label{5.52}
\end{equation}%
\begin{equation}
\left\Vert \nabla h\right\Vert _{C\left( \overline{\Omega }\right) }\leq
C_{1}.  \label{5.53}
\end{equation}

It follows from (\ref{eq:J}) that we need to consider the expression 
\begin{equation}
A=\left\vert L\left( p_{1}+h+F\right) \right\vert ^{2}-\left\vert L\left(
p_{1}+F\right) \right\vert ^{2},  \label{5.54}
\end{equation}%
where the nonlinear operator $L$ is given in (\ref{eq:intdiff}). First, we
will single out the linear, with respect to $h$, part of $A$. This will lead
us to the Frech\'{e}t derivative $J_{\lambda ,\rho }^{\prime }.$ Next, we
will single out $\left\vert \Delta h\right\vert ^{2}.$ This will enable us
to apply Carleman estimate of Theorem \ref{thm:4.1}. \ We have: 
\begin{equation}
\left\vert z_{1}\right\vert ^{2}-\left\vert z_{2}\right\vert ^{2}=\left(
z_{1}-z_{2}\right) \overline{z}_{1}+\left( \overline{z}_{1}-\overline{z}%
_{2}\right) z_{2},\text{ }\forall z_{1},z_{2}\in \mathbb{C}.  \label{5.6}
\end{equation}%
Let%
\begin{equation}
z_{1}=L\left( p_{1}+h+F\right) ,\text{ }z_{2}=L\left( p_{1}+F\right) ,
\label{5.7}
\end{equation}%
Then by (\ref{5.54})-(\ref{5.12})%
\begin{eqnarray}
A_{1} &=&\left( z_{1}-z_{2}\right) \overline{z}_{1},\text{ }A_{2}=\left( 
\overline{z}_{1}-\overline{z}_{2}\right) z_{2},\text{ }  \label{5.8} \\
A &=&A_{1}+A_{2}.  \label{5.81}
\end{eqnarray}%
Taking into account (\ref{eq:intdiff}), (\ref{eq:J}) and (\ref{5.7}), we
obtain%
\begin{equation}
\begin{gathered} z_{1}-z_{2}=\Delta h-2k^{2}\nabla h \left( \nabla
V-\int_{k}^{\overline{k}}\left( \nabla p_{1}+\nabla F\right) d\kappa
\right) \\ +2k\left( \int_{k}^{\overline{k}}\nabla hd\kappa \right)
\left( 2\nabla V-2\int_{k}^{\overline{k}}\left( \nabla p_{1}+\nabla
F\right) d\kappa +k\left( \nabla p_{1}+\nabla F\right) \right) +2i\left(
h_{z}-\int_{k}^{\overline{k}}h_{z}d\kappa \right) . \end{gathered}
\label{5.9}
\end{equation}%
Next,%
\begin{equation*}
\overline{z}_{1}=\left( \Delta \overline{h}+\Delta \overline{p_{1}}+\Delta 
\overline{F}\right) 
\end{equation*}%
\begin{equation*}
-2k\left( \nabla \overline{V}-\int_{k}^{\overline{k}}\left( \nabla 
\overline{p_{1}}+\nabla \overline{h}+\nabla \overline{F}\right) d\kappa
\right) \cdot \left( k\left( \nabla \overline{p_{1}}+\nabla \overline{h}%
+\nabla \overline{F}\right) +\nabla \overline{V}-\int_{k}^{\overline{%
k}}\left( \nabla \overline{p_{1}}+\nabla \overline{h}+\nabla \overline{F}%
\right) d\kappa \right) 
\end{equation*}%
\begin{equation*}
-2i\left( k\left( \overline{p_{1z}}+\overline{h}_{z}+\overline{F_{z}}\right)
+\overline{V_{z}}-\int_{k}^{\overline{k}}\left( \overline{p_{1z}}+%
\overline{h}_{z}+\overline{F_{z}}\right) d\kappa \right) .
\end{equation*}%
Hence, by (\ref{5.8})%
\begin{equation}
A_{1}=\left( z_{1}-z_{2}\right) \overline{z}_{1}=\left\vert \Delta
h\right\vert ^{2}+B_{1}^{\left( linear\right) }\left( h,\mathbf{x},k\right)
+B_{1}\left( h,\mathbf{x},k\right) ,  \label{5.10}
\end{equation}%
where the expression $B_{1}^{\left( linear\right) }\left( h,k\right) $ is
linear with respect to $h=\left( h_{1},h_{2}\right) ,$%
\begin{equation}
\begin{gathered} B_{1}^{\left( linear\right) }\left( h,\mathbf{x},k\right)
=\Delta hG_{1}+\left( \nabla h\nabla G_{2}\right) \cdot G_{3}+\left( \nabla
\overline{h}\nabla G_{4}\right) \cdot G_{5} \\ +G_{7}\cdot \left(
\int_{k}^{\overline{k}}\nabla hd\kappa \right) \nabla
G_{6}+G_{9}\cdot \left( \int_{k}^{\overline{k}}\nabla
\overline{h}d\kappa \right) \nabla G_{8}+G_{10}\left(
h_{z}-\int_{k}^{\overline{k}}h_{z}d\kappa \right) +G_{11}\left(
\overline{h}_{z}-\int_{k}^{\overline{k}}\overline{h}_{z}d\kappa
\right) , \end{gathered}  \label{5.11}
\end{equation}%
where explicit expressions for functions $G_{j}\left( \mathbf{x},k\right)
,j=1,...,11$ can be written in an obvious way. Furthermore, it follows from
these expressions as well as from (\ref{5.52}) that $G_{1},G_{2},G_{4},G_{6}%
\in C_{1}$ and $G_{3},G_{5},G_{7},G_{9},$ $G_{10},G_{11}\in C_{0}$. And also%
\begin{equation}
\left\{ 
\begin{array}{c}
\left\Vert G_{1}\right\Vert _{C_{1}},\left\Vert G_{2}\right\Vert
_{C_{1}},\left\Vert G_{4}\right\Vert _{C_{1}},\left\Vert G_{6}\right\Vert
_{C_{1}}\leq C_{1}, \\ 
\left\Vert G_{3}\right\Vert _{C_{0}},\left\Vert G_{5}\right\Vert
_{C_{0}},\left\Vert G_{7}\right\Vert _{C_{0}},\left\Vert G_{9}\right\Vert
_{C_{0}},\left\Vert G_{10}\right\Vert _{C_{0}},\left\Vert G_{11}\right\Vert
_{C_{0}}\leq C_{1}.%
\end{array}%
\right.   \label{5.12}
\end{equation}%
The term $B_{1}\left( h,k\right) $ in (\ref{5.10}) is nonlinear with respect
to $h$. Applying the Cauchy-Schwarz inequality and also using (\ref{5.52})
and (\ref{5.53}), we obtain%
\begin{equation}
\left\vert B_{1}\left( h,\mathbf{x},k\right) \right\vert \geq \frac{1}{4}%
\left\vert \Delta h\right\vert ^{2}-C_{1}\left\vert \nabla h\right\vert
^{2}-C_{1}\int_{k}^{\overline{k}}\left\vert \nabla h\right\vert
^{2}d\kappa .  \label{5.14}
\end{equation}

Similarly with (\ref{5.10})-(\ref{5.14}) we obtain 
\begin{equation}
A_{2}=\left( \overline{z}_{1}-\overline{z}_{2}\right) z_{2}=B_{2}^{\left(
linear\right) }\left( h,\mathbf{x},k\right) +B_{2}\left( h,\mathbf{x}%
,k\right) ,  \label{5.15}
\end{equation}%
where the term $B_{2}^{\left( linear\right) }\left( h,\mathbf{x},k\right) $
is linear with respect to $h$ and its form is similar with the one of $%
B_{1}^{\left( linear\right) }\left( h,\mathbf{x},k\right) $ in (\ref{5.11}),
although with different functions $G_{j},$ which still satisfy direct
analogs of estimates (\ref{5.12}). As to the term $B_{2}\left( h,\mathbf{x}%
,k\right) ,$ it is nonlinear with respect to $h$ and, as in (\ref{5.14}), 
\begin{equation}
\left\vert B_{2}\left( h,\mathbf{x},k\right) \right\vert \geq \frac{1}{4}%
\left\vert \Delta h\right\vert ^{2}-C_{1}\left\vert \nabla h\right\vert
^{2}-C_{1}\int_{k}^{\overline{k}}\left\vert \nabla h\right\vert
^{2}d\kappa .  \label{5.16}
\end{equation}%
Denote $B\left( h,\mathbf{x},k\right) =B_{1}\left( h,\mathbf{x},k\right)
+B_{2}\left( h,\mathbf{x},k\right) .$ In addition to (\ref{5.14}) and (\ref%
{5.16}), the following upper estimate is valid: 
\begin{equation}
\left\vert B\left( h,\mathbf{x},k\right) \right\vert \leq C_{1}\left(
\left\vert \Delta h\right\vert ^{2}+\left\vert \nabla h\right\vert
^{2}+\int_{k}^{\overline{k}}\left\vert \nabla h\right\vert
^{2}d\kappa \right) .  \label{5.17}
\end{equation}

Thus, it follows from (\ref{eq:intdiff}), (\ref{eq:J}), (\ref{5.7})-(\ref%
{5.81}), (\ref{5.10})-(\ref{5.16}) that 
\begin{equation*}
J_{\lambda ,\rho }\left( p_{1}+h\right) -J_{\lambda ,\rho }\left(
p_{1}\right) =
\end{equation*}%
\begin{equation}
\exp \left( -2\lambda \left( s+d\right) ^{-\nu }\right) \int_{\underline{k}}^{\overline{k}}\int_{\Omega }\left[ S_{1}\Delta
h+S_{2}\cdot \nabla h\right] \varphi _{\lambda }\left( z\right) d\mathbf{x}%
d\kappa +2\rho \left[ h,p_{1}\right]   \label{5.18}
\end{equation}%
\begin{equation*}
+\exp \left( -2\lambda \left( s+d\right) ^{-\nu }\right) \int_{%
\underline{k}}^{\overline{k}}\int_{\Omega }B\left( h,\mathbf{x}%
,\kappa \right) \varphi _{\lambda }\left( z\right) d\mathbf{x}d\kappa .
\end{equation*}%
The second line of (\ref{5.18})%
\begin{equation}
Lin\left( h\right) =\exp \left( -2\lambda \left( s+d\right) ^{-\nu }\right)
\int_{\underline{k}}^{\overline{k}}\int_{\Omega }\left[
S_{1}\Delta h+S_{2}\cdot \nabla h\right] \varphi _{\lambda }\left( z\right) d%
\mathbf{x}d\kappa +2\rho \left[ h,p_{1}\right]   \label{5.19}
\end{equation}%
is linear with respect to $h$, where vector functions $S_{1}\left( \mathbf{x}%
,k\right) ,S_{2}\left( \mathbf{x},k\right) $ are such that 
\begin{equation}
\left\vert S_{1}\left( \mathbf{x},k\right) \right\vert ,\left\vert
S_{2}\left( \mathbf{x},k\right) \right\vert \leq C_{1}\text{ in }\overline{%
\Omega }\times \left[ \underline{k},\overline{k}\right] .  \label{5.20}
\end{equation}%
As to the third line of (\ref{5.18}), it can be estimated from the below as%
\begin{equation}
\begin{gathered} \exp \left( -2\lambda \left( s+d\right) ^{-\nu }\right)
\int_{\underline{k}}^{\overline{k}}\int_{\Omega }B\left(
h,\mathbf{x},\kappa \right) \varphi _{\lambda }\left( z\right)
d\mathbf{x}d\kappa \\ \geq \exp \left( -2\lambda \left( s+d\right) ^{-\nu
}\right) \left[
\frac{1}{2}\int_{\underline{k}}^{\overline{k}}\int_{\Omega
}\left\vert \Delta h\right\vert ^{2}\varphi _{\lambda }\left( z\right)
d\mathbf{x}d\kappa
-C_{1}\int_{\underline{k}}^{\overline{k}}\int_{\Omega
}\left\vert \nabla h\right\vert ^{2}\varphi _{\lambda }\left( z\right)
d\mathbf{x}d\kappa \right] +\rho \left\Vert h\right\Vert _{H_{3}}^{2}.
\end{gathered}  \label{5.21}
\end{equation}%
In addition, (\ref{5.17}) implies that%
\begin{equation}
\begin{gathered} \exp \left( -2\lambda \left( s+d\right) ^{-\nu }\right)
\left\vert \int_{\underline{k}}^{\overline{k}}\int_{\Omega
}B\left( h,\mathbf{x},\kappa \right) \varphi _{\lambda }\left( z\right)
d\mathbf{x}d\kappa \right\vert \\ \leq C_{1}\exp \left( -2\lambda \left(
s+d\right) ^{-\nu }\right)
\int_{\underline{k}}^{\overline{k}}\int_{\Omega }\left(
\left\vert \Delta h\right\vert ^{2}+\left\vert \nabla h\right\vert
^{2}\right) \varphi _{\lambda }\left( z\right) d\mathbf{x}d\kappa +\rho
\left\Vert h\right\Vert _{H_{3}}^{2}. \end{gathered}  \label{5.22}
\end{equation}

First, consider the functional $Lin\left( h\right) $ in (\ref{5.19}). It
follows from (\ref{3.70}), (\ref{5.19}) and (\ref{5.20}) that%
\begin{equation*}
\left\vert Lin\left( h\right) \right\vert \leq C_{1}\exp \left( 2\lambda
t_{\nu }\right) \left\Vert h\right\Vert _{H_{3}}.
\end{equation*}%
Hence, $Lin\left( h\right) :H_{3}\rightarrow \mathbb{R}$ is a bounded linear
functional. Hence, by Riesz theorem for each pair $\lambda ,\nu >0$ there
exists a 2D vector function $Z_{\lambda ,\nu }\in H_{3}$ independent on $h$
such that 
\begin{equation}
Lin\left( h\right) =\left[ Z_{\lambda ,\nu },h\right] ,\text{ }\forall h\in
H_{3}.  \label{5.23}
\end{equation}%
In addition, (\ref{5.17}), (\ref{5.18}) and (\ref{5.23}) imply that%
\begin{equation}
\left\vert J_{\lambda ,\rho }\left( p_{1}+h\right) -J_{\lambda ,\rho }\left(
p_{1}\right) -\left[ Z_{\lambda ,\nu },h\right] \right\vert \leq C_{1}\exp
\left( 2\lambda t_{\nu }\right) \left\Vert h\right\Vert _{H_{3}}^{2}.
\label{5.24}
\end{equation}%
Thus, applying (\ref{5.18})-(\ref{5.24}), we conclude that $Z_{\lambda ,\nu
} $ is the Frech\'{e}t \ derivative of the functional $J_{\lambda ,\rho
}\left( p_{1}\right) $ at the point $p_{1},Z_{\lambda ,\nu }=J_{\lambda
,\rho }^{\prime }\left( p_{1}\right) $.

Thus, (\ref{5.18}) and (\ref{5.21}) imply that 
\begin{equation}  \label{5.25}
\begin{gathered} J_{\lambda ,\rho }\left( p_{1}+h\right) -J_{\lambda ,\rho
}\left( p_{1}\right) -J_{\lambda ,\rho }^{\prime }\left( p_{1}\right) \left(
h\right) \\ \geq \exp \left( -2\lambda \left( s+d\right) ^{-\nu }\right)
\left[
\frac{1}{2}\int_{\underline{k}}^{\overline{k}}\int_{\Omega
}\left\vert \Delta h\right\vert ^{2}\varphi _{\lambda }\left( z\right)
d\mathbf{x}d\kappa
-C_{1}\int_{\underline{k}}^{\overline{k}}\int_{\Omega
}\left\vert \nabla h\right\vert ^{2}\varphi _{\lambda }\left( z\right)
d\mathbf{x}d\kappa \right] +\rho \left\Vert h\right\Vert _{H_{3}}^{2}.
\end{gathered}
\end{equation}%
Assuming that $\lambda \geq \lambda _{0},$ we now apply Carleman estimate of
Theorem \ref{thm:4.1}, 
\begin{equation*}
\frac{1}{2}\int_{\underline{k}}^{\overline{k}}\int_{\Omega
}\left\vert \Delta h\right\vert ^{2}\varphi _{\lambda }\left( z\right) d%
\mathbf{x}d\kappa -C_{1}\int_{\underline{k}}^{\overline{k}%
}\int_{\Omega }\left\vert \nabla h\right\vert ^{2}\varphi _{\lambda
}\left( z\right) d\mathbf{x}d\kappa +\rho \left\Vert h\right\Vert
_{H_{3}}^{2}
\end{equation*}%
\begin{equation*}
\geq \frac{C}{\lambda }\sum_{i,j=1}^{3}\int_{\underline{k}}^{%
\overline{k}}\int_{\Omega }\left\vert h_{x_{i}x_{j}}\right\vert
^{2}\varphi _{\lambda }\left( z\right) d\mathbf{x}d\kappa +C\lambda
\int_{\underline{k}}^{\overline{k}}\int_{\Omega }\left\vert
\nabla h\right\vert ^{2}\varphi _{\lambda }\left( z\right) d\mathbf{x}%
d\kappa -C_{1}\int_{\underline{k}}^{\overline{k}}\int_{%
\Omega }\left\vert \nabla h\right\vert ^{2}\varphi _{\lambda }\left(
z\right) d\mathbf{x}d\kappa +\rho \left\Vert h\right\Vert _{H_{3}}^{2}.
\end{equation*}%
Choosing $\lambda \geq \lambda _{1}$ to be sufficiently large, we obtain%
\begin{equation}  \label{5.26}
\begin{gathered}
\frac{1}{2}\int_{\underline{k}}^{\overline{k}}\int_{\Omega
}\left\vert \Delta h\right\vert ^{2}\varphi _{\lambda }\left( z\right)
d\mathbf{x}d\kappa
-C_{1}\int_{\underline{k}}^{\overline{k}}\int_{\Omega
}\left\vert \nabla h\right\vert ^{2}\varphi _{\lambda }\left( z\right)
d\mathbf{x}d\kappa +\rho \left\Vert h\right\Vert _{H_{3}}^{2}\\ \geq
\frac{C}{\lambda
}\sum_{i,j=1}^{3}\int_{\underline{k}}^{\overline{k}}\int%
_{\Omega }\left\vert h_{x_{i}x_{j}}\right\vert ^{2}\varphi _{\lambda
}\left( z\right) d\mathbf{x}d\kappa +C_{1}\lambda
\int_{\underline{k}}^{\overline{k}}\int_{\Omega }\left\vert
\nabla h\right\vert ^{2}\varphi _{\lambda }\left( z\right)
d\mathbf{x}d\kappa +\rho \left\Vert h\right\Vert _{H_{3}}^{2}. \end{gathered}
\end{equation}%
Finally, noting that $\varphi _{\lambda }\left( z\right) \geq \exp \left(
2\lambda \left( s+d\right) ^{-\nu }\right) $ in $\Omega $ and using (\ref%
{5.25}) and (\ref{5.26}), we obtain%
\begin{equation*}
J_{\lambda ,\rho }\left( p_{1}+h\right) -J_{\lambda ,\rho }\left(
p_{1}\right) -J_{\lambda ,\rho }^{\prime }\left( p_{1}\right) \left(
h\right) \geq \frac{C_{1}}{\lambda }\left\Vert h\right\Vert
_{H_{2}}^{2}+\rho \left\Vert h\right\Vert _{H_{3}}^{2},
\end{equation*}%
\begin{equation*}
J_{\lambda ,\rho }\left( p_{1}+h\right) -J_{\lambda ,\rho }\left(
p_{1}\right) -J_{\lambda ,\rho }^{\prime }\left( p_{1}\right) \left(
h\right) \geq C_{1}\left\Vert h\right\Vert _{H_{1}}^{2}+\rho \left\Vert
h\right\Vert _{H_{3}}^{2}.\text{ \ \ \ \ \ }
\end{equation*}%
$\square $

\subsection{Proof of Theorem \ref{thm:4.4}}

\label{sec:5.3}

This proof is completely similar with the proof of theorem 3.1 of \cite%
{BakushinskiiKlibanov17} and is, therefore, omitted.

\subsection{Proof of Theorem \ref{thm:4.5}}

\label{sec:5.4}

The existence and uniqueness of the minimizer $p_{\min ,\lambda }\in 
\overline{B\left( R\right) }$, inequality (\ref{4.6}) as well as convergence
estimate (\ref{4.90}) follow immediately from the combination of Theorems
\ref{thm:4.3} and \ref{thm:4.4} with lemma 2.1 and theorem 2.1 of \cite{BakushinskiiKlibanov17}. 
$\ \square $

\subsection{Proof of Theorem \ref{thm:4.6}}

\label{sec:5.5}

Temporary denote $L\left( p+F\right) =L\left( p+F,V_{\mu \left( \delta
\right) ,\nu ,\alpha \left( \delta \right) }\right) ,$ $J_{\lambda ,\rho
}\left( p\right) :=$ $J_{\lambda ,\rho }\left( p,F,V_{\mu \left( \delta
\right) ,\nu ,\alpha \left( \delta \right) }\right) $ meaning dependence on
the tail function $V_{\mu \left( \delta \right) ,\nu ,\alpha \left( \delta
\right) }$. Consider the functional $J_{\lambda ,\rho }\left( p,F,V_{\mu
\left( \delta \right) ,\nu ,\alpha \left( \delta \right) }\right) $ for $%
p=p_{\ast },$%
\begin{equation}  \label{5.27}
\begin{gathered} J_{\lambda ,\rho }\left( p_{\ast },F,V_{\mu \left( \delta
\right) ,\nu ,\alpha \left( \delta \right) }\right) =\exp \left( -2\lambda
\left( s+d\right) ^{-\nu }\right)
\int_{\underline{k}}^{\overline{k}}\int_{\Omega }\left\vert
L\left( p_{\ast }+F,V\right) \left( \mathbf{x},\kappa \right) \right\vert
^{2}\varphi _{\lambda }^{2}\left( z\right) d\mathbf{x}d\kappa \\ +\rho
\left\Vert p_{\ast }\right\Vert _{H_{3}}^{2}. \end{gathered}
\end{equation}%
Since $p_{\ast }\in B\left( R\right) $ and $L\left( p_{\ast }+F_{\ast
},V_{\ast }\right) \left( \mathbf{x},\kappa \right) =0,$ then (\ref{5.27})
implies that 
\begin{equation}
J_{\lambda ,\rho }\left( p_{\ast },F_{\ast },V_{\ast }\right) =\rho
\left\Vert p_{\ast }\right\Vert _{H_{3}}^{2}\leq \rho R^{2}=\sqrt{\delta }%
R^{2}.  \label{5.28}
\end{equation}%
It follows from (\ref{eq:intdiff}), (\ref{3.7}), (\ref{3.80}), (\ref{4.100}%
), (\ref{5.27}) and (\ref{5.28}) that 
\begin{equation}
J_{\lambda ,\rho }\left( p_{\ast },F,V_{\mu \left( \delta \right) ,\nu
,\alpha \left( \delta \right) }\right) \leq C_{2}\sqrt{\delta }\sqrt{\ln
\left( \delta ^{-1/\left( 2t_{\nu }\right) }\right) }.  \label{5.29}
\end{equation}%
Next, using (\ref{4.2}) and recalling that $\lambda =\ln \left( \delta
^{-1/\left( 2t_{\nu }\right) }\right) $, we obtain 
\begin{equation*}
\begin{gathered} J_{\lambda ,\rho }\left( p_{\ast },F,V_{\mu \left( \delta
\right) ,\nu ,\alpha \left( \delta \right) }\right) -J_{\lambda ,\rho
}\left( p_{\min ,\lambda \left( \delta \right) },F,V_{\mu \left( \delta
\right) ,\nu ,\alpha \left( \delta \right) }\right) -J_{\lambda ,\rho
}^{\prime }\left( p_{\min ,\lambda \left( \delta \right) },F,V_{\mu \left(
\delta \right) ,\nu ,\alpha \left( \delta \right) }\right) \left( p_{\ast
}-p_{\min ,\lambda }\right) \\ \geq \frac{C_{2}}{\ln \left( \delta
^{-1/\left( 2t_{\nu }\right) }\right) }\left\Vert p_{\ast }-p_{\min ,\lambda
}\right\Vert _{H_{2}}^{2}. \end{gathered}
\end{equation*}%
Next, since $-J_{\lambda ,\rho }\left( p_{\min ,\lambda \left( \delta
\right) },F,V_{\mu \left( \delta \right) ,\nu ,\alpha \left( \delta \right)
}\right) \leq 0$ and also by (\ref{4.6}) $-J_{\lambda ,\rho }^{\prime
}\left( p_{\min ,\lambda },F,V_{\mu \left( \delta \right) ,\nu ,\alpha
\left( \delta \right) }\right) \left( p_{\ast }-p_{\min ,\lambda \left(
\delta \right) }\right) \leq 0,$ we obtain, using (\ref{5.29}),%
\begin{equation*}
\left\Vert p_{\ast }-p_{\min ,\lambda \left( \delta \right) }\right\Vert
_{H_{2}}^{2}\leq C_{2}\sqrt{\delta }\left[ \ln \left( \delta ^{-1/\left(
2t_{\nu }\right) }\right) \right] ^{3/2},
\end{equation*}%
which implies (\ref{4.7}). Estimate (\ref{4.8}) follows immediately from (%
\ref{4.7}), (\ref{3.9}) and Remark 3.1.

We now prove (\ref{4.10}) and (\ref{4.11}). Using (\ref{4.90}), (\ref{4.7})
and the triangle inequality, we obtain for $n=1,2,...$%
\begin{equation*}
\begin{gathered} \left\Vert p_{\ast }-p_{n}\right\Vert _{H_{2}}\leq
\left\Vert p_{\ast }-p_{\min ,\lambda \left( \delta \right) }\right\Vert
_{H_{2}}+\left\Vert p_{\min ,\lambda \left( \delta \right)
}-p_{n}\right\Vert _{H_{2}}\leq C_{2}\delta ^{1/4}\left[ \ln \left( \delta
^{-1/\left( 2t_{\nu }\right) }\right) \right] ^{3/4}+\left\Vert p_{\min
,\lambda \left( \delta \right) }-p_{n}\right\Vert _{H_{3}}\\ \leq
C_{2}\delta ^{1/4}\left[ \ln \left( \delta ^{-1/\left( 2t_{\nu }\right)
}\right) \right] ^{3/4}+r^{n}\left\Vert p_{\min ,\lambda }-p_{0}\right\Vert
_{H_{3}}, \end{gathered}
\end{equation*}%
which proves (\ref{4.10}). Next, using (\ref{4.8}) and (\ref{4.90}), we
obtain 
\begin{equation*}
\begin{gathered} \left\Vert c_{\ast }-c_{n}\right\Vert _{L_{2}\left( \Omega
\right) }\leq \left\Vert c_{\ast }-c_{\min ,\lambda \left( \delta \right)
}\right\Vert _{L_{2}\left( \Omega \right) }+\left\Vert c_{\min ,\lambda
\left( \delta \right) }-c_{n}\right\Vert _{L_{2}\left( \Omega \right) }\\
\leq C_{2}\delta ^{1/4}\left[ \ln \left( \delta ^{-1/\left( 2t_{\nu }\right)
}\right) \right] ^{3/4}+C_{2}\left\Vert p_{\min ,\lambda \left( \delta
\right) }-p_{n}\right\Vert _{H_{2}}\leq C_{2}\delta ^{1/4}\left[ \ln \left(
\delta ^{-1/\left( 2t_{\nu }\right) }\right) \right] ^{3/4}+C_{2}\left\Vert
p_{\min ,\lambda \left( \delta \right) }-p_{n}\right\Vert _{H_{3}}\\ \leq
C_{2}\delta ^{1/4}\left[ \ln \left( \delta ^{-1/\left( 2t_{\nu }\right)
}\right) \right] ^{3/4}+C_{2}r^{n}\left\Vert p_{\min ,\lambda \left( \delta
\right) }-p_{0}\right\Vert _{H_{3}}. \end{gathered}
\end{equation*}%
The latter proves (\ref{4.11}). $\square $

\section{Numerical Study}

\label{sec:6}

In this section, we describe some details of the numerical implementation of
the proposed globally convergent method and demonstrate results of
reconstructions for computationally simulated data. Recall that, as it is
stated in section 2.1, our applied goal in numerical studies is to calculate
locations and dielectric constants of targets which mimic antipersonnel land
mines and IEDs. We model these targets as small sharp inclusions located in
an uniform background, which is air with its dielectric constant $c\left(
air\right) =1.$ Sometimes IEDs can indeed be located in air. In addition, in
previous works \cite{KlibanovLiem17buried,KlibanovThanh15b} of the first
author with coauthors the problem of imaging of targets mimicing land mines
and IEDs in the case when those targets are buried in a sandbox is
considering. Microwave experimental data are used in these publications. The
tail functions numerical method was used in these works. It was demonstrated
in \cite{KlibanovLiem17buried,KlibanovThanh15b} that, after applying certain
data preprocessing procedures, one can treat those targets as ones located
in air. Recall that $c\left( air\right) =1$ is a good approximation for the
value of the dielectric constant of air. Thus, in this paper, we conduct
numerical experiments for the case when small inclusions of our interest are
located in air. We test several of values of the dielectric constant and
sizes of those inclusions. However, we do not assume in computations the
knowledge of the background in the domain of interest $\Omega $ in (\ref%
{eq:2.1}), except of the knowledge that $c\left( \mathbf{x}\right) =1$
outside of $\Omega ,$ see (\ref{eq:coef}).

\subsection{The Carleman Weight Function of numerical studies}

\label{sec:6.1}

The CWF $\varphi _{\lambda }\left( z\right) =\exp \left( 2\lambda \left(
z+s\right) ^{-\nu }\right) ,$ which was introduced in (\ref{3.3}), changes
too rapidly due to the presence of the parameter $\nu >0.$ We have
established in our computational experiments that such a rapid change does
not allow us to obtain good numerical results, also, see page 1581 of \cite%
{Baud} for a similar conclusion. Hence, we use in our numerical studies a
simpler CWF $\psi _{\lambda }\left( z\right) ,$%
\begin{equation}
\psi _{\lambda }\left( z\right) =e^{-2\lambda z}.  \label{6.1}
\end{equation}%
We cannot prove an analog of Theorem \ref{thm:4.1} for this CWF. Nevertheless, the
following Carleman estimate is valid in the 1D case \cite{KlibanovKolesov17}:%
\begin{equation}
\int_{-\xi }^{d}\left( w^{\prime \prime }\right) ^{2}\psi _{\lambda
}\left( z\right) dz\geq C_{3}\left[ \int_{-\xi }^{d}\left( w^{\prime
\prime }\right) ^{2}\psi _{\lambda }\left( z\right) dz+\lambda
\int_{-\xi }^{d}\left( w^{\prime }\right) ^{2}\psi _{\lambda }\left(
z\right) dz+\lambda ^{3}\int_{-\xi }^{d}w^{2}\psi _{\lambda }\left(
z\right) dz\right] ,  \label{6.2}
\end{equation}%
for for all $\lambda >1$ and for any real valued function $w\in H^{2}\left(
-\xi ,d\right) $ such that $w\left( -\xi \right) =w^{\prime }\left( -\xi
\right) =0.$ Here and below the number $C_{3}=C_{3}\left( \xi ,d\right) >0$
depends only on numbers $\xi $ and $d$.

To briefly justify the CWF (\ref{6.1}) from the analytical standpoint,
consider now the case when the Laplace operator is written in partial finite
differences with respect to variables $x,y\in \left[ -b,b\right] $ (see (\ref%
{eq:2.1})) with the uniform grid step size $h>0$ with respect to each
variable $x$ and $y,$ 
\begin{equation}
\Delta ^{h}=\frac{\partial ^{2}}{\partial z^{2}}+\Delta _{x,y}^{h}.
\label{6.20}
\end{equation}%
Here $\Delta _{x,y}^{h}$ is the Laplace operator with respect to $x,y$,
which is written in finite differences. Suppose that we have $M_{h}$
interior grid points in each direction $x$ and $y$. The domain $\Omega $ in (%
\ref{eq:2.1}) becomes 
\begin{equation*}
\Omega _{h}=\left\{ \left( x_{j},y_{s},z\right) :\left\vert x_{j}\right\vert
,\left\vert y_{s}\right\vert <b,z\in \left( -\xi ,d\right) \right\} ;\text{ }%
j,s=1,...,M_{h},
\end{equation*}%
where $\left( x_{j},y_{s}\right) $ are grid points. Then the finite
difference analog of the integral of $\left( \Delta u\right) ^{2}\psi
_{\lambda }\left( z\right) $ over the domain $\Omega $ is%
\begin{equation}
Z_{h}\left( u,\lambda \right)
=\sum_{j,s=1}^{M_{h}}h^{2}\int_{-\xi }^{d}\left[ \left(
u_{zz}+u_{xx}^{h}+u_{yy}^{h}\right) \left( x_{j},y_{s},z\right) \right]
^{2}\psi _{\lambda }\left( z\right) dz,  \label{6.3}
\end{equation}%
where $u\left( x_{j},y_{s},z\right) $ is the discrete real valued function
defined in $\Omega _{h}$ and such that $u_{zz}\left( x_{j},y_{s},z\right)
\in L_{2}\left( -\xi ,d\right) $ for all $\left( x_{j},y_{s}\right) .$ In
addition, $u\left( -\xi ,x_{j},y_{s}\right) =\partial _{z}u\left( -\xi
,x_{j},y_{s}\right) =0.$ Also, in (\ref{6.3}) $u_{xx}^{h}$ and $u_{yy}^{h}$
are corresponding finite difference derivatives of the function $u\left(
x_{j},y_{s},z\right) $ at the point $\left( x_{j},y_{s},z\right) $.
\textquotedblleft Interior" grid points are those located in $\overline{%
\Omega }\diagdown \partial \Omega .$ As to the grid points located at $%
\partial \Omega ,$ they are counted in the well known way in finite
differences derivatives in (\ref{6.3}). Obviously,%
\begin{equation}
Z_{h}\left( u,\lambda \right) \geq \frac{1}{2}\sum%
_{j,s=1}^{M_{h}}h^{2}\int_{-\xi }^{d}\left[ u_{zz}\left(
x_{j},y_{s},z\right) \right] ^{2}\psi _{\lambda }\left( z\right) dz-\widehat{%
C}\sum_{j,s=1}^{M_{h}}\int_{-\xi }^{d}\left[ u\left(
x_{j},y_{s},z\right) \right] ^{2}\psi _{\lambda }\left( z\right) dz.
\label{6.4}
\end{equation}%
Here and below in this section the number $\widehat{C}=\widehat{C}\left(
1/h\right) >0$ depends only on $1/h.$ Hence, the following analog of the
Carleman estimate (\ref{6.2}) for the case of the operator (\ref{6.20})
follows immediately from (\ref{6.4}): 
\begin{equation}  \label{6.5}
\begin{gathered} Z_{h}\left( u,\lambda \right) \geq
C_{3}\sum_{j,s=1}^{M_{h}}h^{2}\int_{-\xi }^{d}\left[
u_{zz}\left( x_{j},y_{s},z\right) \right] ^{2}\psi _{\lambda }\left(
z\right) dz\\ +\widehat{C}\left[ \lambda
\sum_{j,s=1}^{M_{h}}\int_{-\xi }^{d}\left[ u_{z}\left(
x_{j},y_{s},z\right) \right] ^{2}\psi _{\lambda }\left( z\right) dz+\lambda
^{3}\sum_{j,s=1}^{M_{h}}\int_{-\xi }^{d}\left[ u\left(
x_{j},y_{s},z\right) \right] ^{2}\psi _{\lambda }\left( z\right) dz\right]
,\forall \lambda \geq \widetilde{\lambda }\left( h\right) >1, \end{gathered}
\end{equation}%
where $\widetilde{\lambda }\left( h\right) $ increases with the decrease of $%
h$.

Suppose now that operators $\Delta $ and $\nabla $ in (\ref{eq:intdiff}), (%
\ref{eq:Jtail}), (\ref{eq:J}) are rewritten in partial finite differences
with respect to $x,y.$ As to the spaces $H^{3}(\Omega )$ and $H_{3+r},$ they
were introduced to ensure that functions $p\in C_{1},V\in C^{1}\left( 
\overline{\Omega }\right) ,F\in C_{2},$ see (\ref{3.73}), (\ref{3.74}). Note
that by the embedding theorem $H^{n}\left( -\xi ,d\right) \subset C^{n-1}%
\left[ -\xi ,d\right] ,n\geq 1$. Thus, we replace the space $H^{m}(\Omega )$
with $m=1,2,3$ in (\ref{3.700}) with the following finite difference analog
of it for complex valued functions $f$:%
\begin{equation*}
H^{n,h}(\Omega _{h})=\left\{ f\left( x_{j},y_{s},z\right) :\left\Vert
f\right\Vert _{H^{n,h}(\Omega
_{h})}^{2}=\sum_{j,s=1}^{M_{h}}\sum_{r=0}^{n}h^{2}\int%
_{-\xi }^{d}\left\vert \partial _{z}^{r}f\left( x_{j},y_{s},z\right)
\right\vert ^{2}dz\right\} ,n=1,2,
\end{equation*}%
and similarly for the replacement of $H_{m}$ with $H_{n,h}.$ So, we replace
the regularization terms $\alpha \Vert V\Vert _{H^{3}(\Omega )}^{2}$ and $%
\rho \left\Vert p\right\Vert _{H_{3}}^{2}$ in (\ref{eq:Jtail}) and (\ref%
{eq:J}) with $\alpha \Vert V_{h}\Vert _{H^{2,h}(\Omega )}^{2}$ and $\rho
\left\Vert p_{h}\right\Vert _{H_{2,h}}^{2}$ respectively. Also, we replace
in (\ref{3.80}) $H_{4}$ with $H_{3,h}$ and in (\ref{3.10}) we replace $H_{3}$
with $H_{2,h}.$ The functionals $J_{\lambda ,\rho }\left( p+F\right) $ and $%
\widetilde{I}_{\mu ,\alpha }\left( W\right) $ in (\ref{eq:J}) and (\ref{5.1}%
) are replaced with their finite difference analogs, 
\begin{equation}
\widetilde{I}_{\mu ,\alpha }^{h}\left( W_{h}\right) =\exp \left( 2\mu
d\right) \int_{\Omega }\left\vert \Delta ^{h}W_{h}+\Delta
^{h}Q_{h}\right\vert ^{2}\psi _{\mu }\left( z\right) d\mathbf{x}+\alpha
\Vert W_{h}+Q_{h}\Vert _{H^{2,h}(\Omega _{h})}^{2},\text{ }  \label{6.50}
\end{equation}%
\begin{equation}
J_{\lambda ,\rho }^{h}\left( p_{h}\right) =\exp \left( 2\lambda d\right)
\int_{\underline{k}}^{\overline{k}}\int_{\Omega }\left\vert
L^{h}\left( p_{h}+F_{h}\right) \left( \mathbf{x},\kappa \right) \right\vert
^{2}\varphi _{\lambda }^{2}\left( z\right) d\mathbf{x}d\kappa +\rho
\left\Vert p_{h}\right\Vert _{H_{2,h}}^{2},  \label{6.51}
\end{equation}%
where $V_{h},W_{h},p_{h},Q_{h}$ and $F_{h}$ are finite difference analogs of
functions $V,W,p,Q$ and $F$ respectively and $L^{h}$ is the finite
difference analog of the operator $L$ in which operators $\Delta $ and $%
\nabla $ are replaced with their above mentioned finite difference analogs.

Then the Carleman estimate (\ref{6.5}) implies that the straightforward
analogs of Theorems \ref{thm:4.2}-\ref{thm:4.6} are valid for functionals (\ref{6.50}) and (\ref%
{6.51}).\ The only restriction is that the grid step size $h$ should be
bounded from the below, 
\begin{equation}
h\geq h_{0}=const.>0.  \label{6.6}
\end{equation}%
In other words, numerical experiments should not be conducted for the case
when $h$ tends to zero, as it is done sometimes for forward problems for
PDEs. It is our computational experience that condition (\ref{6.6}) is
sufficient for computations. So, we do not change $h$ in our numerical
studies below.

\textbf{Remarks 6.1: }

\begin{enumerate}
\item For brevity, we do not reformulate here those analogs of Theorems
\ref{thm:4.2}-\ref{thm:4.6}. Also, both for brevity and convenience we describe our procedures
below for the case of the continuous spatial variable $\mathbf{x}$. Still,
we actually work in our computations with functionals (\ref{6.50}) and (\ref%
{6.51}).

\item The reason why we have presented the above theory for the case of the
CWF (\ref{3.3}) is that it is both consistent and is valid for the 3D case.
We believe that this theory is interesting in its own right from the
analytical standpoint. On the other hand, in the case of the CWF (\ref{6.1})
and the assumption about partial finite differences, the corresponding
theory (unlike computations!) is similar to the one which we (with
coauthors) have developed in the 1D case of \cite{KlibanovKolesov17}.
\end{enumerate}

\subsection{Data simulation and propagation}

\label{sec:6.2}

To computationally simulate the boundary data $g_{0}(\mathbf{x},k)$ in (\ref%
{eq:cisp}), we solve the Lippmann-Schwinger integral equation 
\begin{equation}
u(\mathbf{x},k)=e^{ikz}+k^{2}\int_{\Omega }\Phi (\mathbf{x},\mathbf{y},k)(c(\mathbf{y} )-1)u(\mathbf{y}k)d\mathbf{y},  \label{eq:LippmannSchwinger}
\end{equation}%
where $\Phi (\mathbf{x},\mathbf{y},k)$ is the fundamental solution of the
Helmholtz equation with $c(\mathbf{x})\equiv 1$: 
\begin{equation*}
\Phi (\mathbf{x},\mathbf{y},k)=\frac{e^{ik|\mathbf{x}-\mathbf{y}|}}{4\pi |\mathbf{x}-%
\mathbf{y}|},\quad \mathbf{x}\neq \mathbf{y}.
\end{equation*}%
The spectral method of \cite{Vainikko00}, which is based on the
periodization technique and the fast Fourier transform, is used to solve (%
\ref{eq:LippmannSchwinger}), see, e.g. \cite{LechleiterLiem14} for the
numerical implementation of this method in MATLAB.

We work with dimensionless variables. Typically linear sizes of
antipersonnel land mines are between 5 and 10 centimeters (cm), see, e.g. 
\cite{Landmine}. Hence, just as in papers with experimental data of our
research group \cite{KlibanovLiem17buried,KlibanovLiem17exp}, we introduce
the dimensionless variables $\mathbf{x}^{\prime }=\mathbf{x}/(10\,\text{cm})$%
. Our mine-like targets are ball-shaped.\ Hence, their radii $r=0.3$ and $%
0.5,$ for example, correspond to diameters of those balls of 6 cm and 10 cm
respectively. This change of variables leads to the dimensionless frequency $%
k$, which is also called the \textquotedblleft wavenumber". Hereafter, for
convenience and brevity, \ we will leave the same notations for
dimensionless spatial variables $\mathbf{x}$ as before. Note that the
dimensionless wavenumber $k=16.2,$ which we work with below (see (\ref{6.7}%
)), corresponds to the frequency of $f=7.7$ GHz. Since in \cite%
{KlibanovLiem17buried,KlibanovLiem17exp} microwave experimental data were
collected by our research group for the range of frequencies from 1 GHz to
10 GHz, then $f=7.7$ GHz is a realistic value of the frequency.

\begin{table}[tbp]
\caption{Mine-like and IED-like inclusions tested in our numerical studies.
A single inclusion in cases 1-3. Two inclusions simultaneously: the left inclusion is 4.1,
the right inclusion is 4.2}
\label{tab1}
\begin{center}
\begin{tabular}{ccc}
\hline
Inclusion number & $\max \left( c\right) $ in the inclusion & Radius $r$ \\ 
\hline
1 & 3 & 0.3 \\ 
2 & 3 & 0.5 \\ 
3 & 5 & 0.3 \\ 
4.1 & 7 & 0.3 \\ 
4.2 & 3 & 0.5 \\ \hline
\end{tabular}%
\end{center}
\end{table}

All inclusions, which we have numerically tested, are listed in Table \ref%
{tab1}, where $r$ denotes the radius of the corresponding ball-shaped
inclusion. To have smooth target/background interfaces, the dielectric
constants of inclusions were smoothed out for a better stability of the
numerical method of solving the Lippmann-Schwinger equation (\ref%
{eq:LippmannSchwinger}). But the maximal values of dielectric constants
remain unchanged in this smoothing, and these values are reached in the
centers of those balls. In our study, the center of each ball representing a
single inclusion is at the point $\mathbf{x}=\left( x,y,z\right) =\left(
0,0,0\right) $ and centers of two inclusions in case number 4 of Table \ref%
{tab1} are placed at points $\left( x,y,z\right) =(-0.75,0,0)$ (left
inclusion) and $\left( x,y,z\right) =(0.75,0,0)$ (right inclusion). However,
when running the inversion procedure, we do not assume the knowledge of
neither those centers nor the shapes of those inclusions.

\begin{figure}[tbp]
\begin{center}
\subfloat[\label{fig:f_noiseless}]{\includegraphics[width=0.33%
\textwidth]{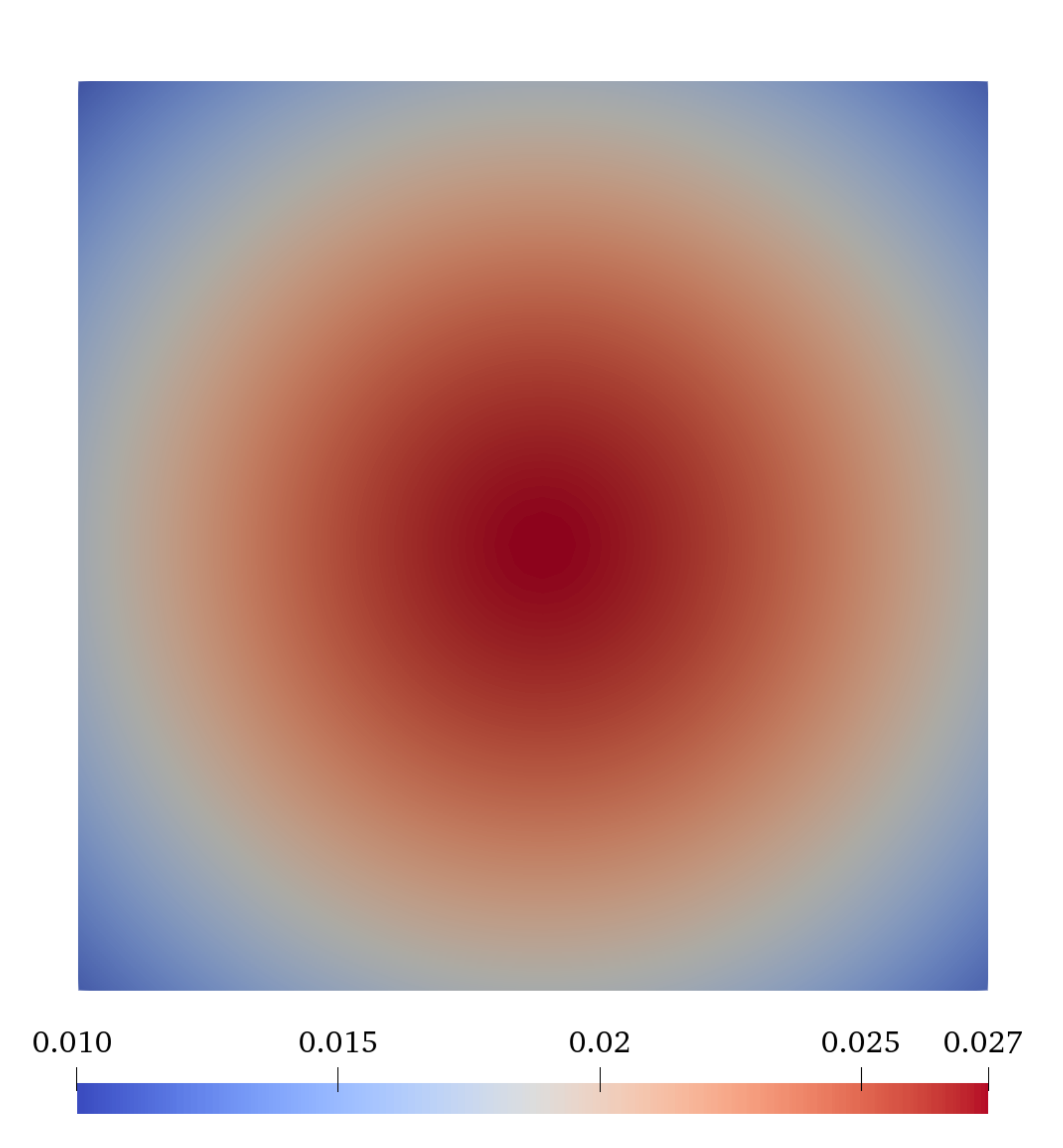}} 
\subfloat[\label{fig:g_noiseless}]{%
\includegraphics[width=0.33\textwidth]{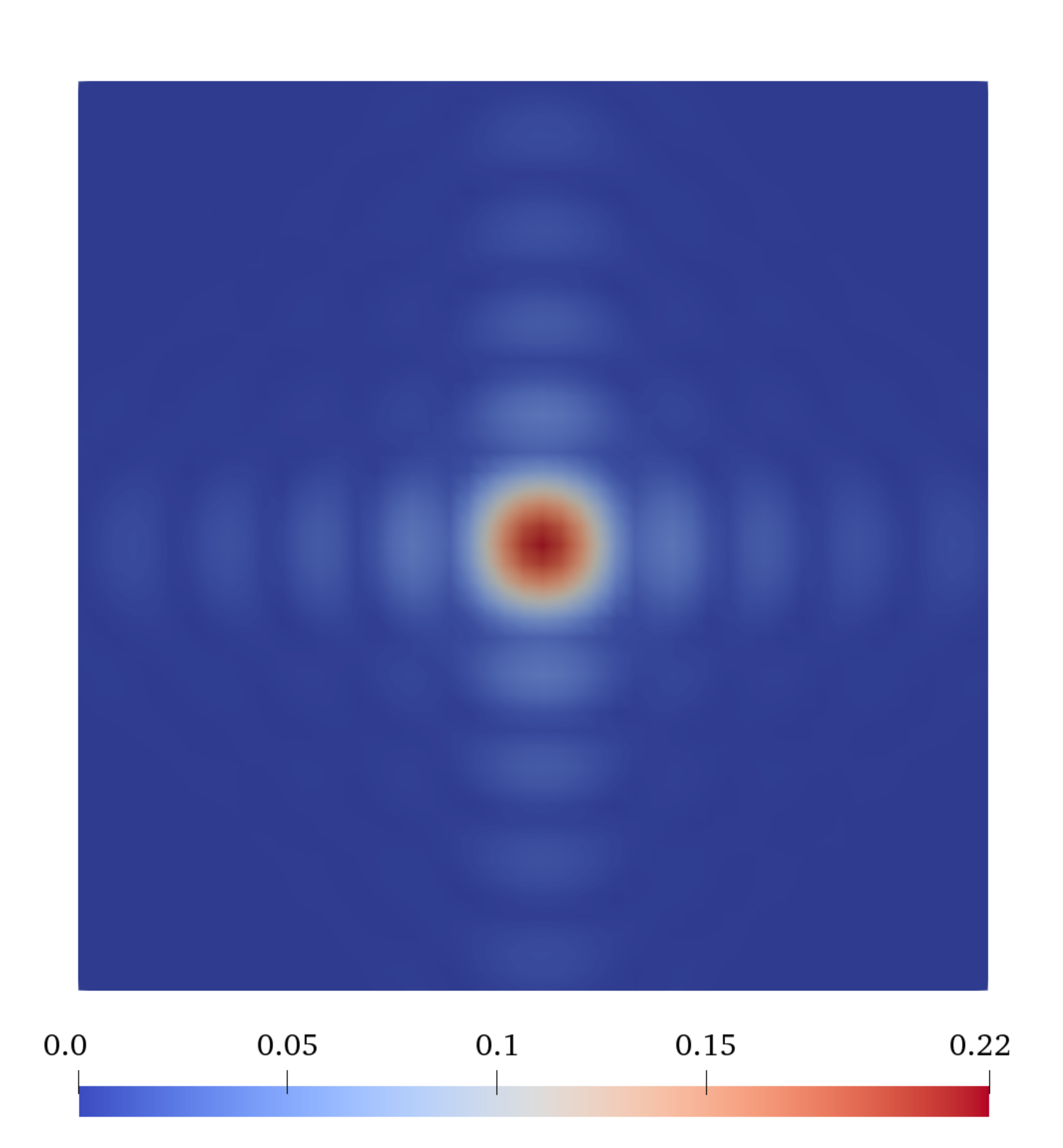}} 
\subfloat[%
\label{fig:q_noiseless}]{\includegraphics[width=0.33\textwidth]{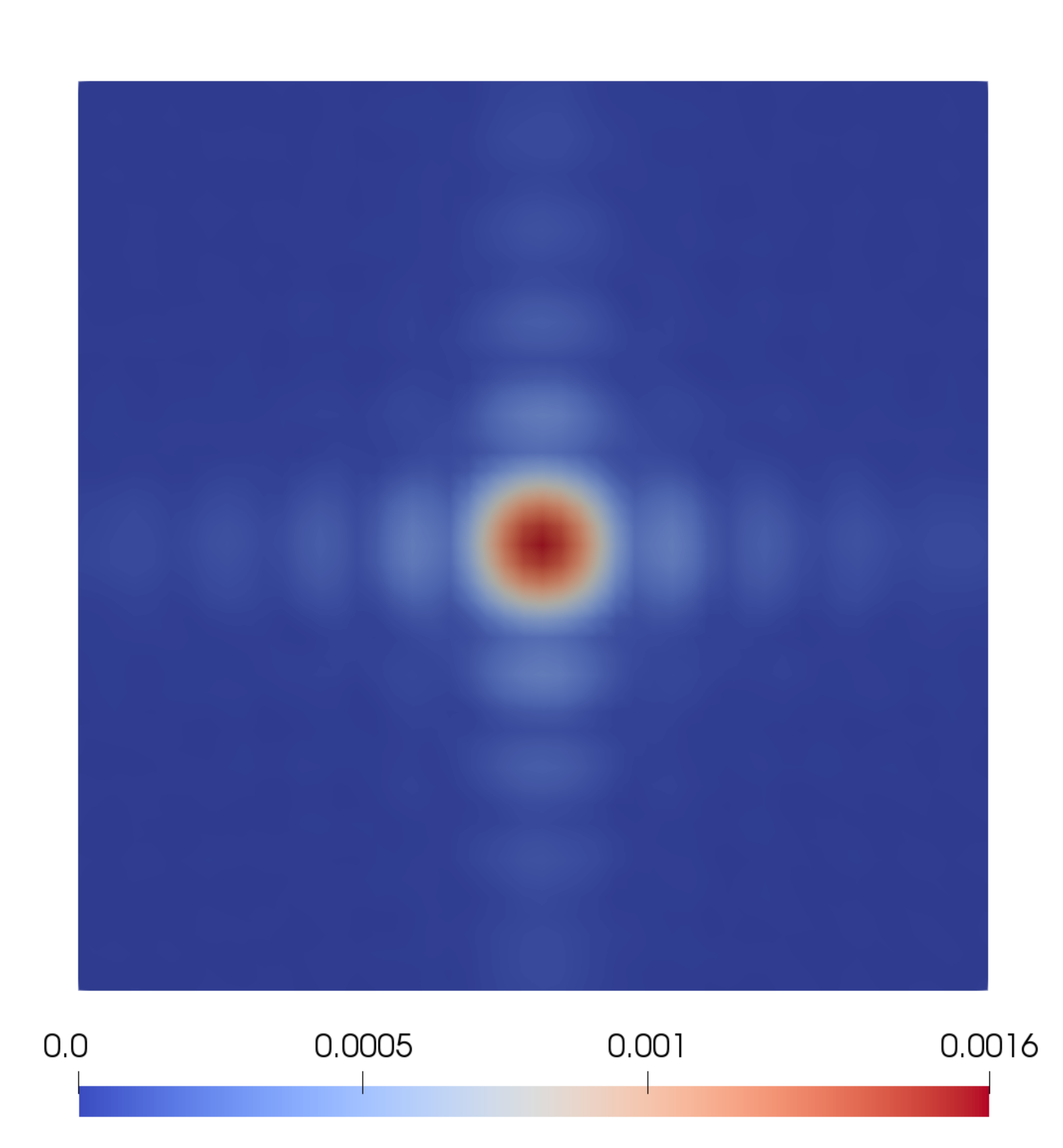}}
\end{center}
\caption{The absolute value of the simulated data without the noise for the
inclusion 1 in Table \protect\ref{tab1}. Here $k = 16.2$. see (\protect\ref%
{6.7}). a) The measured data $f(\mathbf{x}, k)$, b) the propagated data $g_0(%
\mathbf{x}, q)$, c) the function $\protect\phi_0(\mathbf{x},k)$. Hence, $x,y$
coordinates of inclusion are clearly seen on the propagated data, unlike the
measured data. This indicates that the data propagation procedure is quite a
useful one.}
\label{fig:noiseless_data}
\end{figure}

In the setup for our computational experiments, we want to be close to the
experimental setup of \cite{KlibanovLiem17buried,KlibanovLiem17exp}.
Actually, in \cite{KlibanovLiem17buried,KlibanovLiem17exp} the data are
collected not at the part $\Gamma $ (\ref{eq:2.1}) of the boundary of the
domain $\Omega $ as in (\ref{eq:cisp}). Instead, they are collected on a
square $P_{meas}$, which is a part of the so-called measurement plane 
\begin{equation}
P_{m}=\{z=-A\},  \label{60}
\end{equation}
where $A=const.>\xi .$ We solve the Lippmann-Schwinger equation (\ref%
{eq:LippmannSchwinger}) to obtain computationally simulated data $f(\mathbf{x%
},k)$ for $\mathbf{x}\in P_{meas}.$ We refer to $f(\mathbf{x},k)$ as the
\textquotedblleft measured data". The measurement plane $P_{m}$ is located
far from $\Gamma $. This causes several complications. First, we would need
to solve our CISP in a large computational domain, which could be
time-consuming process. Second, looking at the measured data is not clear
enough how to distinguish inclusions; see Fig. \ref{fig:f_noiseless}.

Hence, we need to propagate the measured data $f(\mathbf{x},k)$ generated by
the Lippmann-Schwinger solver from the rectangle $P_{meas}$ to the so-called
propagation plane $P_{p}=\{z=A^{\prime }\}$, $A^{\prime }\leq -\xi ,$ which
is closer to our inclusions. In fact, we propagate to the plane which
includes the rectangle $\Gamma .$ As a result we get the so-called
propagated data (Fig. \ref{fig:g_noiseless}), which are more focused on the
target of our interest than the original data. So, we can clearly see now
the location of our inclusion in $x,y$ coordinates. The resulting function $%
u\left( \mathbf{x},k\right) $ is our given boundary data $g_{0}(\mathbf{x},k)
$ in (\ref{eq:cisp}) for our CISP. The derivative $u_{z}\left( \mathbf{x}%
,k\right) $ for $\mathbf{x}\in \Gamma ,$ i.e. the function $g_{1}(\mathbf{x}%
,k)$ in (\ref{eq:gz0}), is calculated by propagating $f(\mathbf{x},k)$ into
a plane $\left\{ z=-\xi -\varepsilon \right\} $ for a small number $%
\varepsilon >0$. Next, the finite difference is used to approximate $g_{1}(%
\mathbf{x},k)$.

For brevity, we do not describe the data propagation procedure here.
Instead, we refer to \cite{KlibanovKolesov17exp, KlibanovLiem17exp,
KlibanovLiem17buried} for detailed descriptions. In fact, this procedure is
quite popular in Optics under the name the \emph{angular spectrum
representation method} \cite{Nov}.

We also remark that both in the data propagation procedure and in our
convexification method we need to calculate the some derivatives of noisy
data: the $\partial _{k,z}^{2}-$derivative of the propagated data and the $%
\partial _{k}-$derivative in the convexification. In all cases this is done
using finite differences. We have not observed any instabilities probably
because the step sizes of our finite differences were not too small. The
same was in all previous above cited publications of this research group.

\begin{figure}[tbp]
\begin{center}
\includegraphics[width=0.7\textwidth]{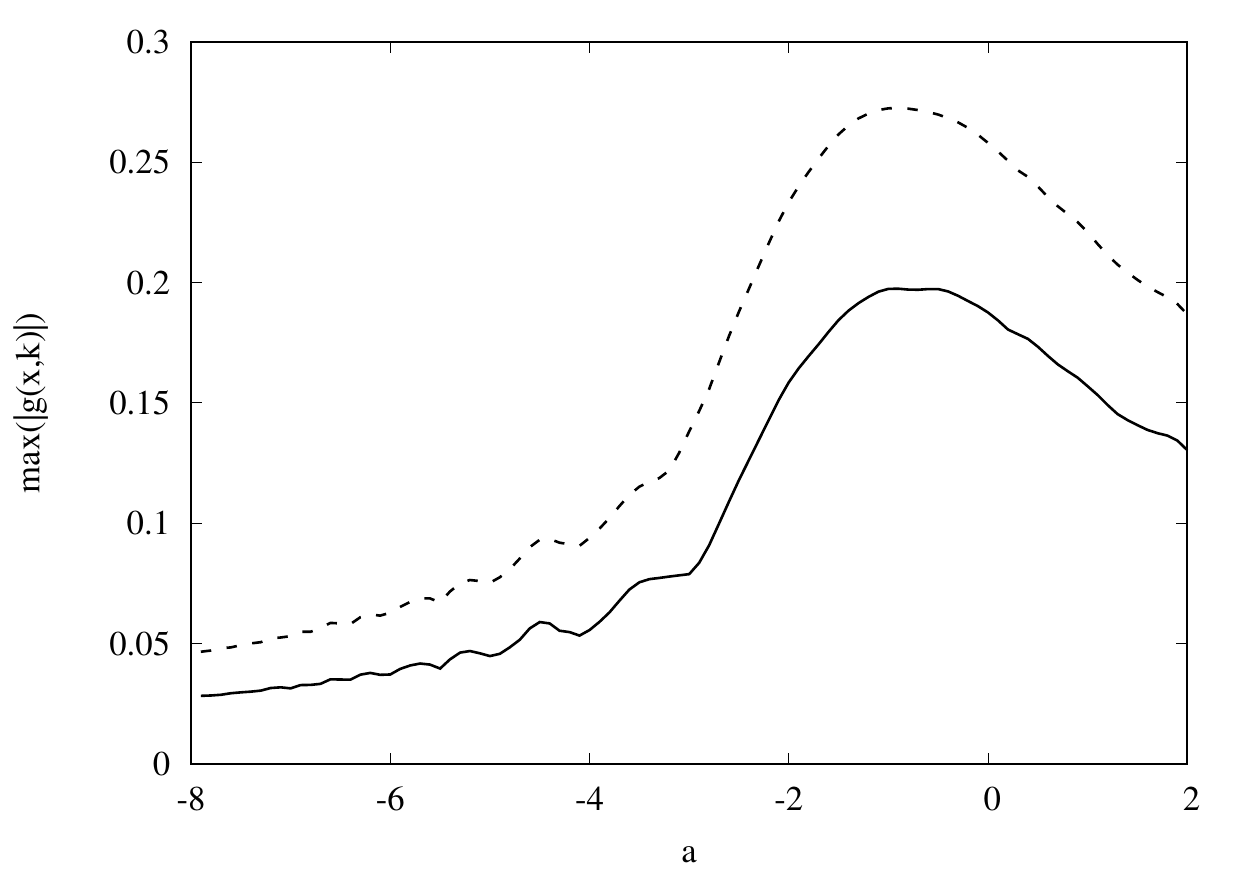}
\end{center}
\caption{The depedence of maximum absolute value of propagated data $g(\mathbf{x}, k)$ on the locaton of propagation plane $a$ for inclusions with $c=3.0$
(solid line) and $c = 5.0$ (dashed line).}
\label{fig:g_max}
\end{figure}

To propagate the function $f(\mathbf{x},k)$ close to an inclusion, we need
to figure out first where this inclusion is located, i.e. we need to
estimate the number $\xi $ in (\ref{eq:2.1}). Fortunately, the data
propagation procedure allows us to do this. For example, consider two
inclusions with the same size $r=0.3$, but with different dielectric
constants $c=3$ and $c=5 $. Centers of both are located at the point $\left(
0,0,0\right) .$ We solve the Lippmann-Schwinger equation for each of these
two cases to generate the data at the measurement plane $P_{m}$ with $A=8,$
see (\ref{60}). Next, we propagate the data to several propagated planes $%
P_{p,a}=\left\{ z=a\right\} ,$ where $a=\left( -8,2\right] $. Here, we use $%
k=16.2$ (see (\ref{6.7})). 
The dependence of the maximal absolute value of the propagated data $M\left(
a\right) =\max_{P_{p,a}}\left\vert u\left( x,y,a,16.2\right) \right\vert $
on the number $a$ for these inclusions is depicted in Fig. \ref{fig:g_max}.
We see that the function $M\left( a\right) $ attains its maximal value near
the point $a_{0}=-0.5$ for both cases. This point is located reasonably
close to the actual position of the front faces (at $z=-0.15)$ of the
corresponding inclusions. The function $M\left( a\right) $ has attained its
maximal value at points $a$ close to $a_{0}$ for all other inclusions we
have tested. Therefore, we propagate the measured data for all inclusions to
the propagated plane $P_{p,-0.5}=\left\{ z=-0.5\right\} $ and we set in (\ref%
{eq:2.1}) $\xi =-0.5$.

We have found in our computations that the optimal interval of wavenumbers
is: 
\begin{equation}
k\in \left[ 15.2,16.2\right] .  \label{6.7}
\end{equation}%
We divide this interval in ten (10) subintervals with the step size $\Delta
k=0.1.$ For each $k=15.2,15.3,...16.1,16.2$ and for each inclusion under the
consideration\ we solve Lippmann-Schwinger equation (\ref%
{eq:LippmannSchwinger}) to generate the function $f(\mathbf{x},k).$ Next, by
propagating this data, we obtain the functions $g_{0}(\mathbf{x},k)$ and $%
g_{1}(\mathbf{x},k)$ in (\ref{eq:cisp}) and (\ref{eq:gz0}) respectively.
Using (\ref{eq:uinc}), (\ref{eq:w}), (\ref{eq:v}) and (\ref{eq:q}), consider
the function $q(\mathbf{x},k)$ on the propagated plane $P_{p}$, i.e. at the
boundary $\Gamma $. In fact, this function is denoted as $\phi _{0}\left( 
\mathbf{x},k\right) $ in (\ref{eq:intdiffbcs}) and it is one of the two
boundary conditions (the second one is $\phi _{1}\left( \mathbf{x},k\right) $
in (\ref{eq:intdiffbcs})) which generate the function $F_{h}$ in the
functional $J_{\lambda ,\rho }^{h}\left( p_{h}\right) $ in (\ref{6.51}).
Fig. \ref{fig:q_noiseless} displays the function $\phi _{0}(\mathbf{x},k)$
for the inclusion number 1 in Table 1 for $k=16.2$.


\subsection{Computational domain}

\label{sec:6.3}

To model the experimental setup of \cite%
{KlibanovLiem17buried,KlibanovLiem17exp} we use the following measurement
plane: 
\begin{equation*}
P_{m}=\left\{ z=-8\right\} ,P_{meas}=\{\mathbf{x}:(x,y)\in (-3,3)\times
(-3,3),z=-8\},
\end{equation*}%
where $P_{meas}\subset P_{m}$ is the square on which measurements are
conducted and $z=-8$ corresponds to the 80 cm. The latter is the approximate
distance from the center of any inclusion to the plane $\left\{ z=0\right\} $
where detectors are located in \cite{KlibanovLiem17buried,KlibanovLiem17exp}%
. Solving equation (\ref{eq:LippmannSchwinger}), we generate the function $f(%
\mathbf{x},k),k\in \left[ 15.2,16.2\right] .$ Next, we propagate this
function to $\Gamma ,$ 
\begin{equation*}
\Gamma =\{\mathbf{x}:(x,y)\in (-3,3)\times (-3,3),z=-0.5\}\subset P_{p}.
\end{equation*}%
Here, $z=-0.5$ was found in section 6.2. Finally, we define our
computational domain as 
\begin{equation}
\Omega =\{\mathbf{x}:(x,y,z)\in (-3,3)\times (-3,3)\times (-0.5,4.5)\}.
\label{eq:omega}
\end{equation}

\subsection{Adding noise}

\label{sec:6.4}

\begin{figure}[tbp]
\begin{center}
\subfloat[\label{fig:f_noisy}]{\includegraphics[width=0.33%
\textwidth]{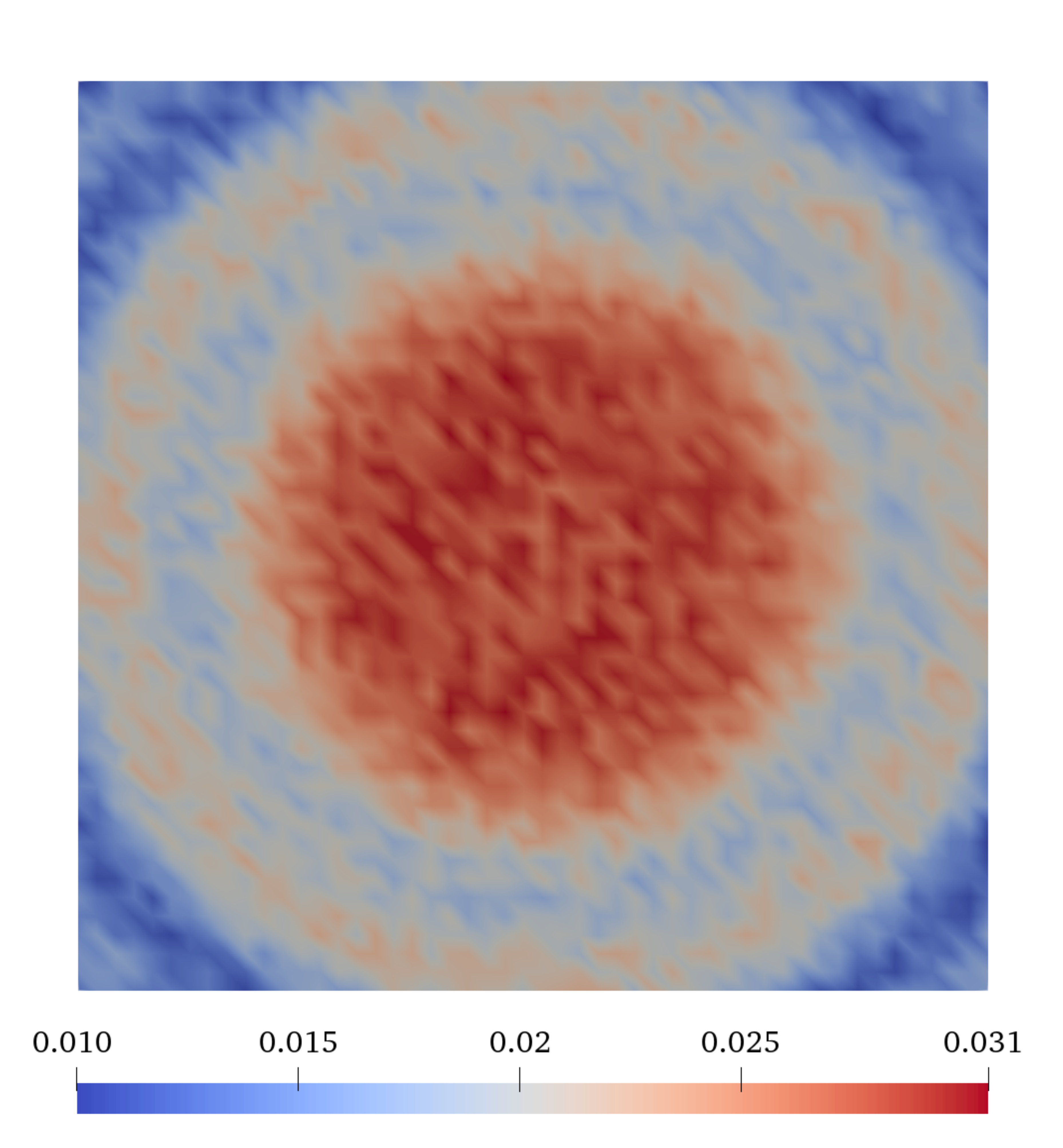}} 
\subfloat[\label{fig:g_noisy}]{%
\includegraphics[width=0.33\textwidth]{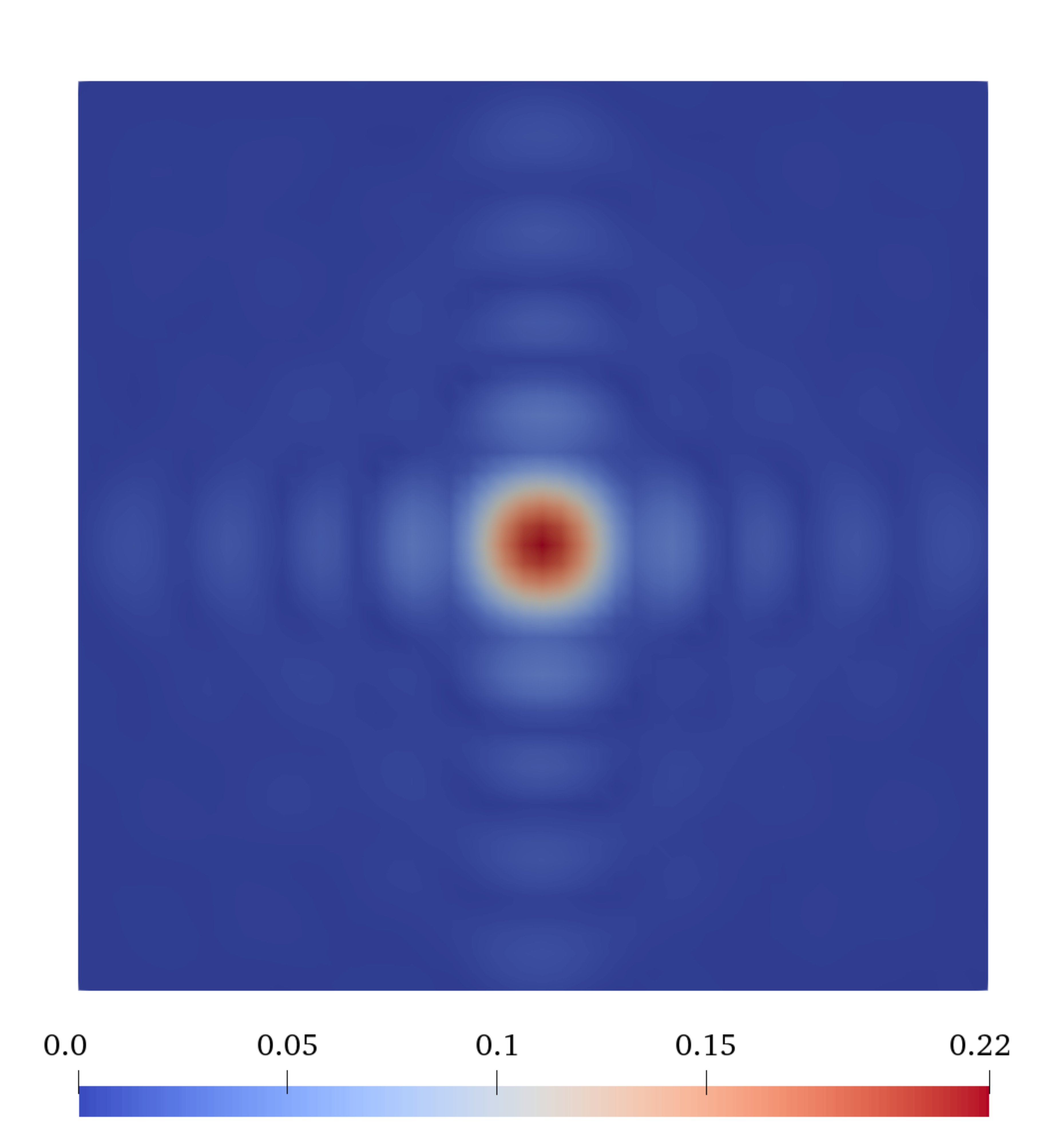}} 
\subfloat[%
\label{fig:q_noisy}]{\includegraphics[width=0.33\textwidth]{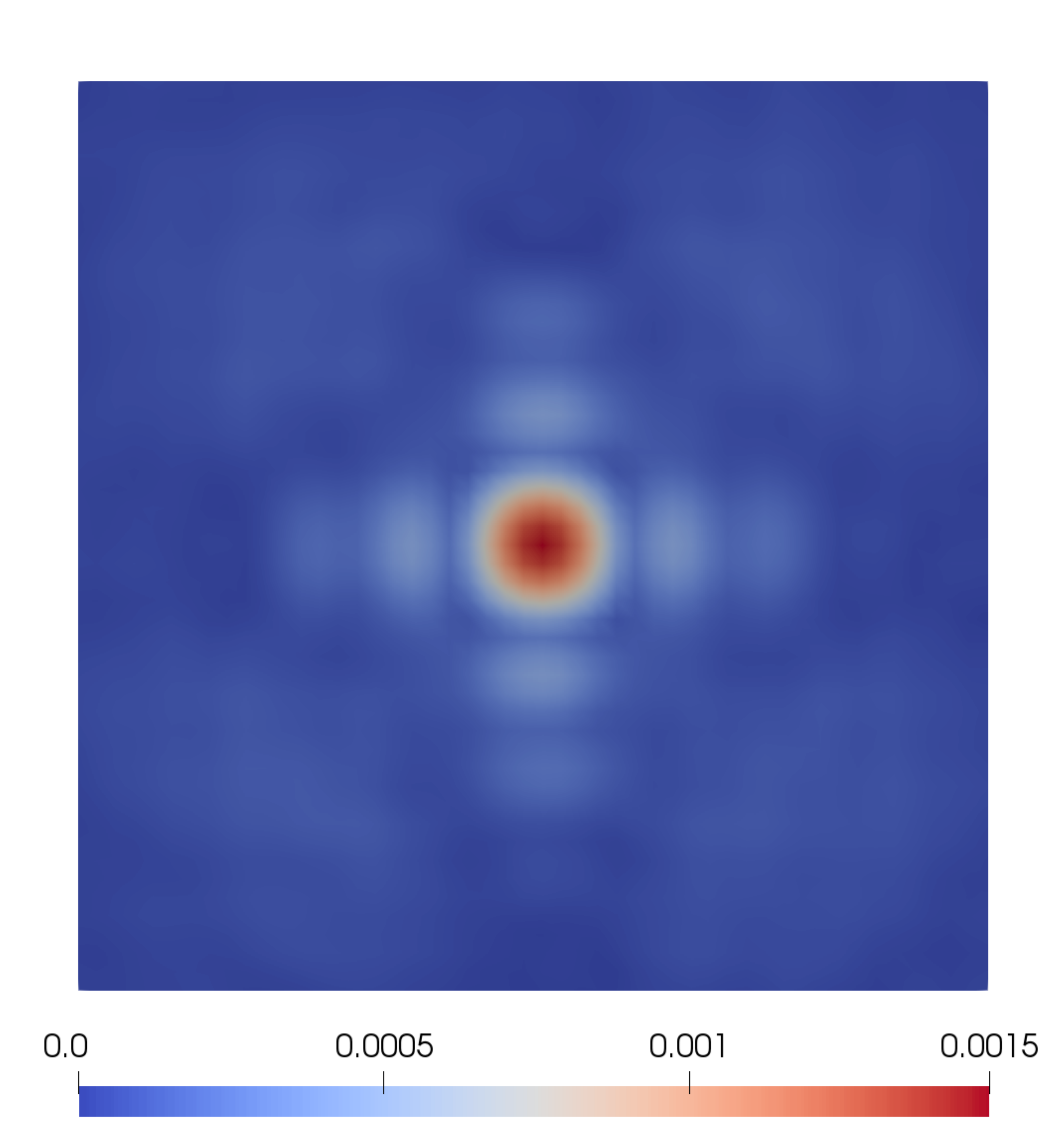}}
\end{center}
\caption{The absolute value of the computationally simulated data with noise
for the inclusion 1 in Table \protect\ref{tab1}. Here $k=16.2$, see (\protect
\ref{6.7}). a) The measured data $f_{noisy}(\mathbf{x},k)$, b) the
propagated data $g_{0,noisy}(\mathbf{x},k)$, c) and the function $\protect\phi %
_{0,noisy}(\mathbf{x},k)$). Comparison of b) with Fig. \protect\ref%
{fig:g_noiseless} shows that the data propagation procedure has a
significant smoothing effect. Also, just as in Fig. \protect\ref%
{fig:noiseless_data}, comparison of b) and c) with a) shows that this
procedure helps quite significantly to find $x,y$ coordinates of targets of
interest.}
\label{fig:noisy_data}
\end{figure}

We add a random noise to the simulated data $f(\mathbf{x},k)$ as follows: 
\begin{equation*}
f_{noisy}(\mathbf{x},k)=f(\mathbf{x},k)+\delta \Vert f(\mathbf{x},k)\Vert
_{L^{2}(\Gamma )}\frac{\sigma (\mathbf{x},k)}{\Vert \sigma (\mathbf{x}%
,k)\Vert _{L^{2}(\Gamma )}}.
\end{equation*}%
Here, $\delta $ is the noise level. Next, $\sigma (\mathbf{x},k)=\sigma _{1}(%
\mathbf{x},k)+i\sigma _{2}(\mathbf{x},k)$, where $\sigma _{1}(\mathbf{x},k)$ and $%
\sigma _{2}(\mathbf{x},k)$ are random numbers uniformly distributed on the
interval $(-1,1)$. We use below $\delta =0.15$, i.e. $15\%$ of the additive
noise. Fig. \ref{fig:noisy_data} displays the absolute value of simulated
data with noise $f_{noise}(\mathbf{x},k)$, the corresponding propagated data 
$g_{0,noisy}(\mathbf{x},k)$ and the function $\phi _{0,noisy}(\mathbf{x},k)$
for the same inclusion and the wavenumber $k=16.2$ as in Fig. \ref%
{fig:noiseless_data}. We see that the data propagation procedure has a
smoothing effect on our noisy measured data, since $g_{0}(\mathbf{x},k)$ in
Fig. \ref{fig:g_noiseless} and $g_{0,noisy}(\mathbf{x},k)$ in Fig. \ref%
{fig:g_noisy} are almost identical.

\subsection{The algorithm}

\label{sec:6.5}

Based on the above theory, we use the following algorithm for determining
the function $c(\mathbf{x})$ from simulated data with noise $f(\mathbf{x},k)$
(here, the subscript \textquotedblleft $noisy$" is left out for convenience,
also see item 1 in Remarks 6.1):

\begin{enumerate}
\item Using the data propagation procedure, calculate the boundary data $%
g_{0}(\mathbf{x},k)$ and $g_{1}(\mathbf{x},k)$.

\item Calculate the subsequent boundary conditions $\phi _{0}(\mathbf{x},k)$%
, $\phi _{1}(\mathbf{x},k)$, $\psi _{0}(\mathbf{x})$, and $\psi _{1}(\mathbf{%
x})$.

\item Compute the auxiliary functions $Q_{h}(\mathbf{x})$ and $F_{h}(\mathbf{%
x},k)$.

\item Compute the minimizer $W_{\min ,h}(\mathbf{x})$ of the functional $%
\widetilde{I}_{\mu ,\alpha }^{h}\left( W_{h}\right) $ in (\ref{6.50}).

\item Using the computed function $V_{h}(\mathbf{x})=W_{\min ,h}(\mathbf{x}%
)+Q_{h}(\mathbf{x})$, minimize the functional $J_{\lambda ,\rho }^{h}(p_{h})$ in (%
\ref{6.51}). Let the function $p_{h,\min }(\mathbf{x},k)$ be its minimizer.
Calculate the function $q_{h}(\mathbf{x},k)=p_{h,\min }(\mathbf{x},k)+F_{h}(%
\mathbf{x},k)$.

\item Compute the function $v_{h}(\mathbf{x},k)$ for $k=\underline{k}$ as
follows: 
\begin{equation*}
v_{h}(\mathbf{x},\underline{k})=-\int_{\underline{k}}^{\overline{k}}q_{h}(%
\mathbf{x},\kappa )d\kappa +V_{h}(\mathbf{x}).
\end{equation*}

\item Calculate the approximation for the unknown coefficient $c(\mathbf{x})$
using the following formulae, see (\ref{eq:coef}), (\ref{eq:intdiffv}) 
\begin{equation*}
\beta (\mathbf{x})=-\Delta ^{h}v_{h}(\mathbf{x},\underline{k})-\underline{k}%
^{2}\nabla v_{h}(\mathbf{x},\underline{k})\cdot \nabla v(\mathbf{x},%
\underline{k})+2i\underline{k}v_{z}(\mathbf{x},\underline{k}),
\end{equation*}%
\begin{equation*}
c\left( \mathbf{x}\right) =\left\{ 
\begin{array}{c}
\mathop{\rm Re}\beta \left( \mathbf{x}\right) +1,\text{ if }\mathop{\rm Re}%
\beta \left( \mathbf{x}\right) \geq 0\text{ and }\mathbf{x}\in \Omega , \\ 
1,\text{ otherwise.}%
\end{array}%
\right.
\end{equation*}
\end{enumerate}

\subsection{Numerical implementation}

\label{sec:6.6}

We now present some details of the numerical implementation. When minimizing
functionals $\widetilde{I}_{\mu ,\alpha }^{h}\left( W_{h}\right) $ and $%
J_{\lambda ,\rho }^{h}\left( p_{h}\right) $ in (\ref{6.50}) and (\ref{6.51}%
), we use finite differences not only with respect to $x,y$ but with respect
to $z$ as well. Thus, $z-$derivatives in these functionals are also written
in finite differences.\ For brevity we use the same notations $\widetilde{I}%
_{\mu ,\alpha }^{h}\left( W_{h}\right) $ and $J_{\lambda ,\rho }^{h}\left(
p_{h}\right) $ for these functionals.

This is the fully discrete case, unlike the semi-discrete case of (\ref{6.50}%
), (\ref{6.51}). The theory for the fully discrete cases of nonlinear
ill-posed problems for PDEs is not yet developed well.\ It seems that such a
theory is much more complicated than the one for the semi-discrete case.
There are only a few results for the fully discrete case, and all are for
linear ill-posed problems for PDEs, as opposed to our nonlinear case, see,
e.g. \cite{Burman,KS}. Since it is not yet clear to us how to extend above
theorems for the fully discrete case, we are not concerned with such
extensions here.

We minimize resulting functionals with respect to the values of
corresponding functions at grid points. In the computational domain (\ref%
{eq:omega}), we use the uniform grid with $N_{x}=N_{y}=N_{z}=51$ points with
the corresponding step sizes $h_{x},h_{y},h_{z}$, where $h_{x}=h_{y}=h.$ The
grid point labeled $(j,s,l)$ corresponds to $\mathbf{x}=(x,y,z)=(x_{j},y_{s},z_{l})$. In addition, the interval $k=[\underline{k},%
\overline{k}]$ of wavenumbers is divided into $N_{k}=11$ points $k_{n}$ with
the step size $h_{k}$. Hence, we use the following discrete functions $W_{h}(%
\mathbf{x})=W(x_{j},y_{s},z_{l})=W_{j,s,l}$ and $p_{h}(\mathbf{x},k)=p_{h}(x_{j},y_{s},z_{l},k_{n})=p_{j,s,l,n}$ at grid points.

To minimize the functionals $\widetilde{I}_{\mu ,\alpha }^{h}(W_{j,s,l})$
and $J_{\lambda ,\rho }^{h}(p_{j,s,l,n}),$ we use the conjugate gradient
method (CG) instead of the gradient projection method, which is suggested by
our theory. Indeed, similarly with \cite{KlibanovKolesov17}, we have
observed that the results obtained by both these methods are practically the
same. On the other hand, CG is easier to implement numerically than the
gradient projection method. Note that we do not employ the standard line
search algorithm for determining the step size of the CG. Instead, we start
with the step size $10^{-4}$, which is reduced two times if the value of the
corresponding functional on the current iteration exceeds its value on the
previous iteration otherwise it remains the same. The minimization algorithm
is stopped when the step size is less then $10^{-10}$. We use zero as the
starting point of the CG for both functions $W_{j,s,l}$ and $p_{j,s,l,n}$.

Gradients of both functionals $\widetilde{I}_{\mu ,\alpha }^{h}(W_{j,s,l})$
and $J_{\lambda ,\rho }^{h}(p_{j,s,l,n})$ are calculated analytically on
each step, and we do not provide details of this for brevity. Rather, we
refer to formulae (7.7) and (7.8) of \cite{KlibanovKuzhuget10}, where
gradients of similar functionals are calculated analytically using the
Kronecker delta function. Also, due to the difficulty with the numerical
implementation of the $H^{2,h}(\Omega _{h})-$norm, we use the simpler $L_{2}$
norm in (\ref{6.50}). As to (\ref{6.51}), we have established numerically
that the minimization of the functional $J_{\lambda ,\rho }^{h}\left(
p_{h}\right) $ works better if the regularization term is absent. Hence, we
set $\rho =0$ in (\ref{6.51}).

\subsection{Reconstruction results}

\label{sec:6.7}

In this section we present the results of our reconstructions for the
inclusions listed in Table \ref{tab1} using the above algorithm. These
results are obtained using the Carleman Weight Function (\ref{6.1}) with $%
\mu =8$ in (\ref{6.50}) and $\lambda =8$ in (\ref{6.51}). We have found that
these are optimal values of the parameters $\mu $ and $\lambda .$ Table \ref%
{tab3} lists each inclusion with the maximal value $c_{exact}$ of the exact
coefficient $c_{exact}=\max_{inclusion}c\left( \mathbf{x}\right) $, radius $%
r $, the maximal value of the computed coefficient $c_{comp}=%
\max_{inclusion}c(\mathbf{x})$, the relative computational error 
\begin{equation*}
\varepsilon =\frac{|c_{comp}-c_{exact}|}{c_{exact}}\cdot 100\%,
\end{equation*}%
and location, i.e. the $z$ coordinate of the point where the value of $%
c_{comp}$ is achieved.

\begin{figure}[tbp]
\begin{center}
\subfloat[\label{fig:c_true3}]{\includegraphics[width=0.5%
\textwidth]{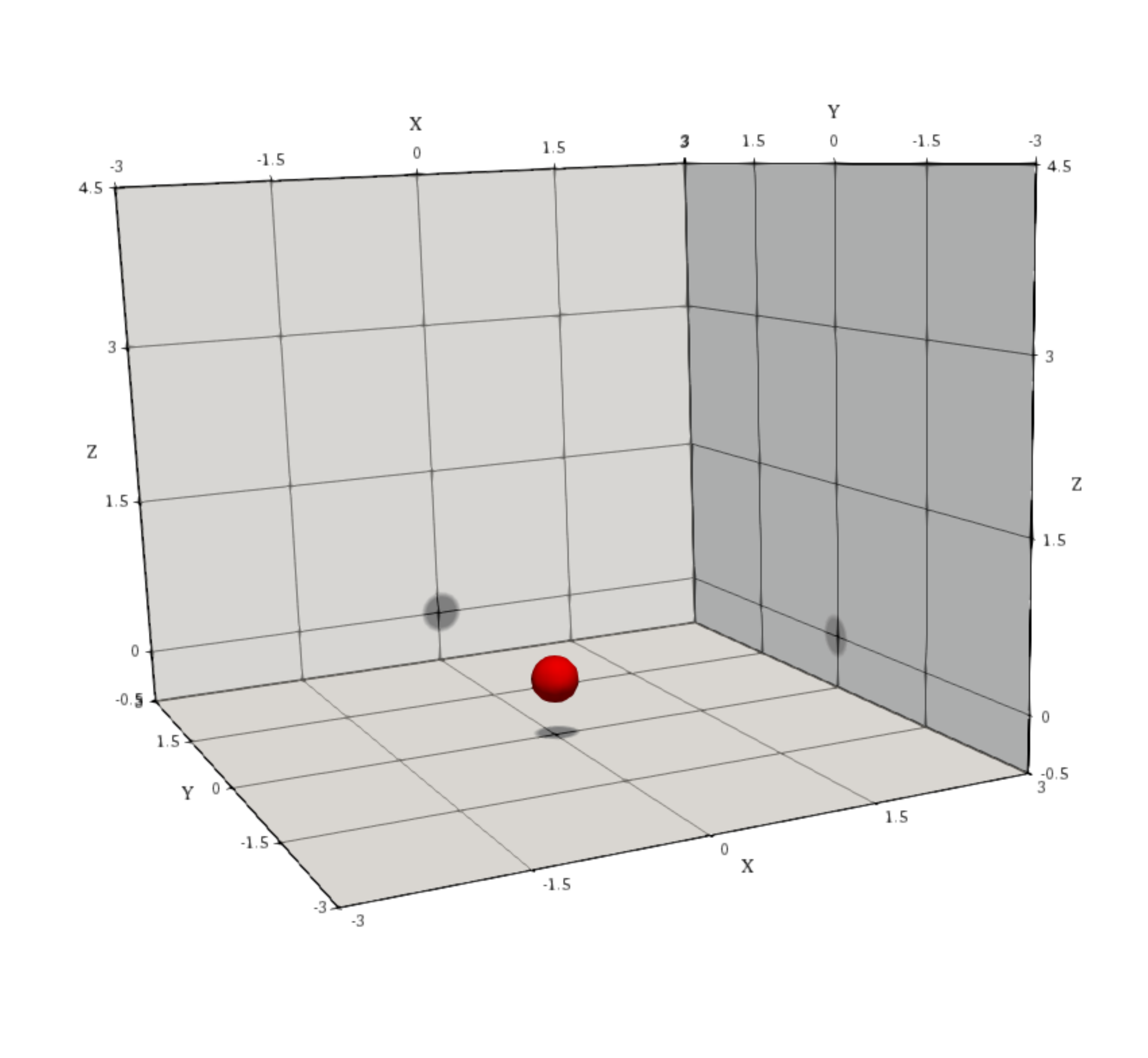}} 
\subfloat[\label{fig:c_comp3}]{%
\includegraphics[width=0.5\textwidth]{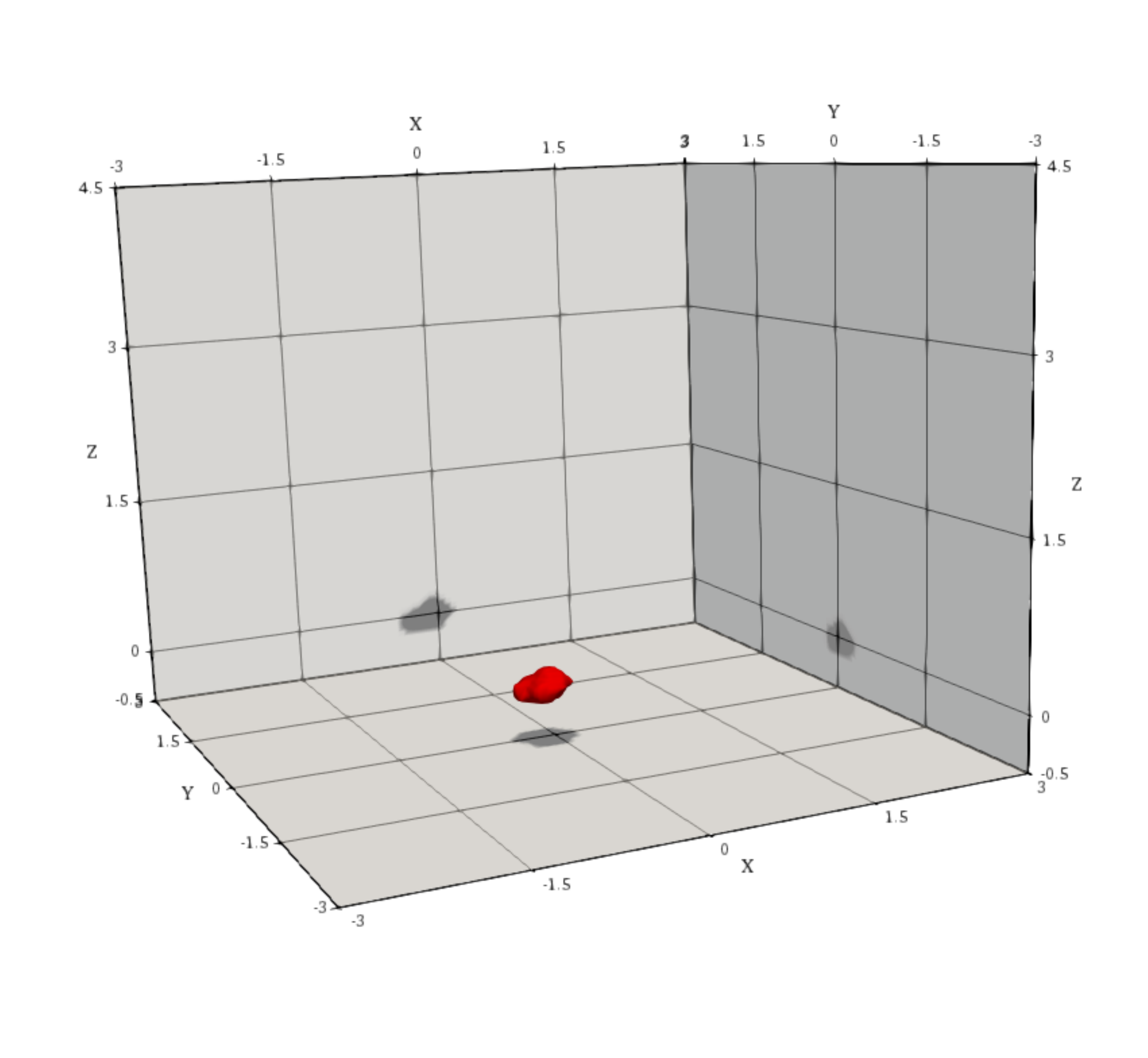}}
\end{center}
\caption{Reconstruction result for the inclusion number 1: exact image (a)
and computed image (b)}
\label{fig:c_3}
\end{figure}

Note that while we have added $15\%$ noise in our simulated data, the
relative computational errors of reconstructed coefficients do not exceed
9\% in all cases, which is 1.67 times less than the level of noise in the
data. Moreover, the locations of points where the values of $c_{comp}$ are
achieved, are reconstructed with a good accuracy as well. Indeed, we need
our reconstructed inclusions to be somewhere between $-r$ and $r$, where
either $r=0.3$ or $r=0.5$. Fig. \ref{fig:c_3} displays the exact and
computed images for the inclusion number 1 in Table \ref{tab1}. Images are
obtained in Paraview.

\begin{figure}[tbp]
\begin{center}
\subfloat[\label{fig:f_2inc}]{\includegraphics[width=0.33\textwidth]{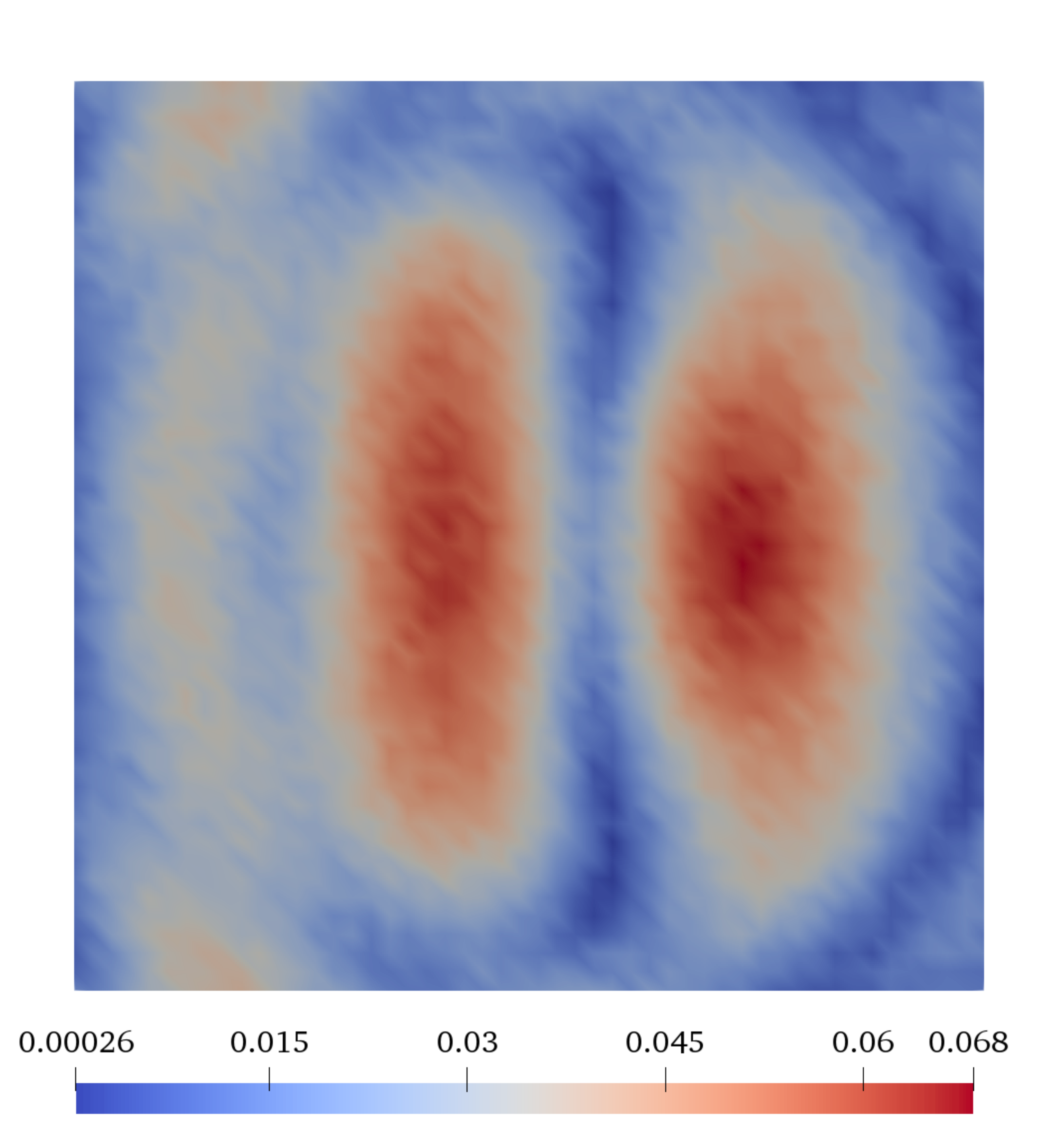}}
\subfloat[\label{fig:g_2inc}]{\includegraphics[width=0.33\textwidth]{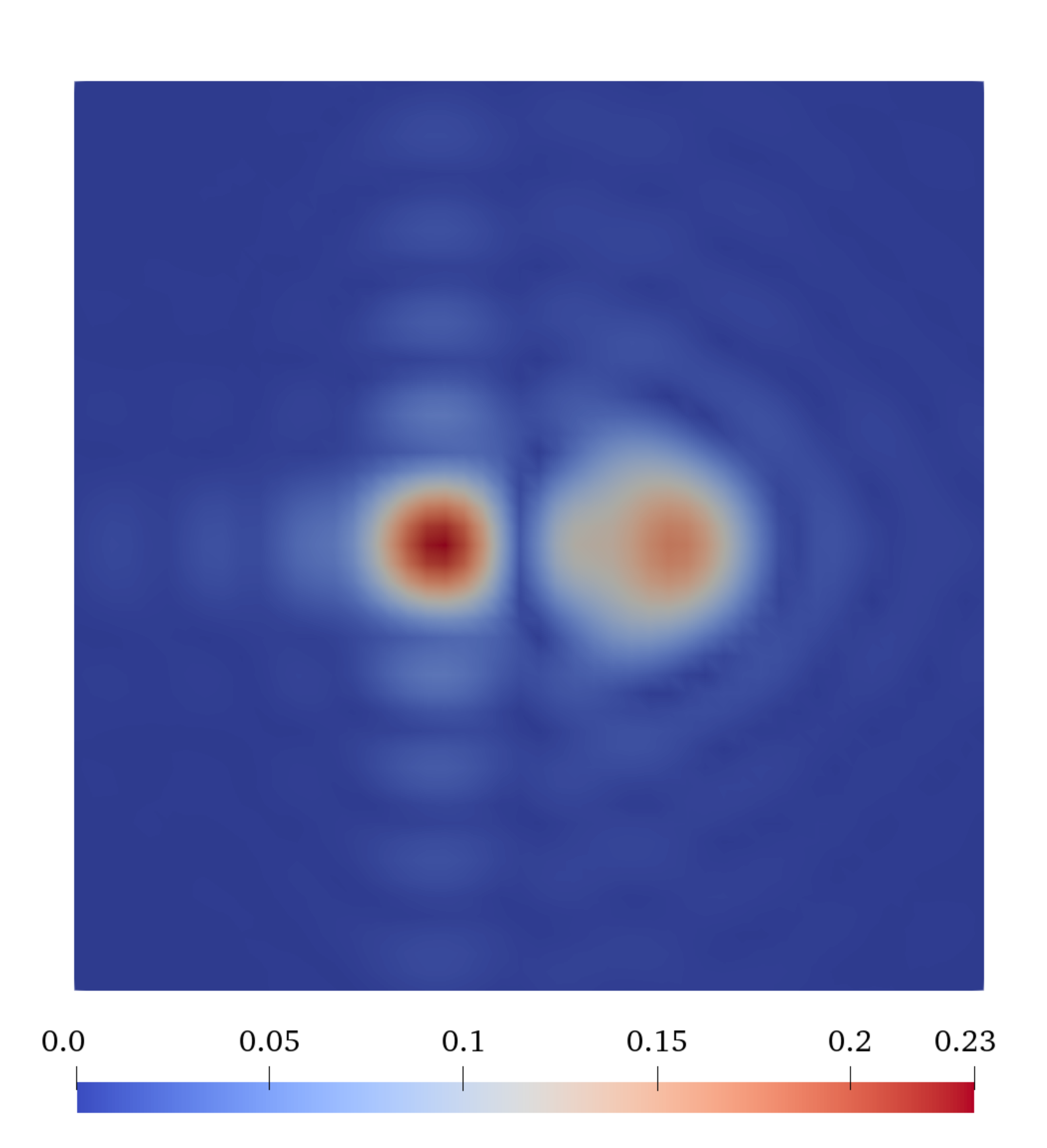}}
\subfloat[\label{fig:q_2inc}]{\includegraphics[width=0.33\textwidth]{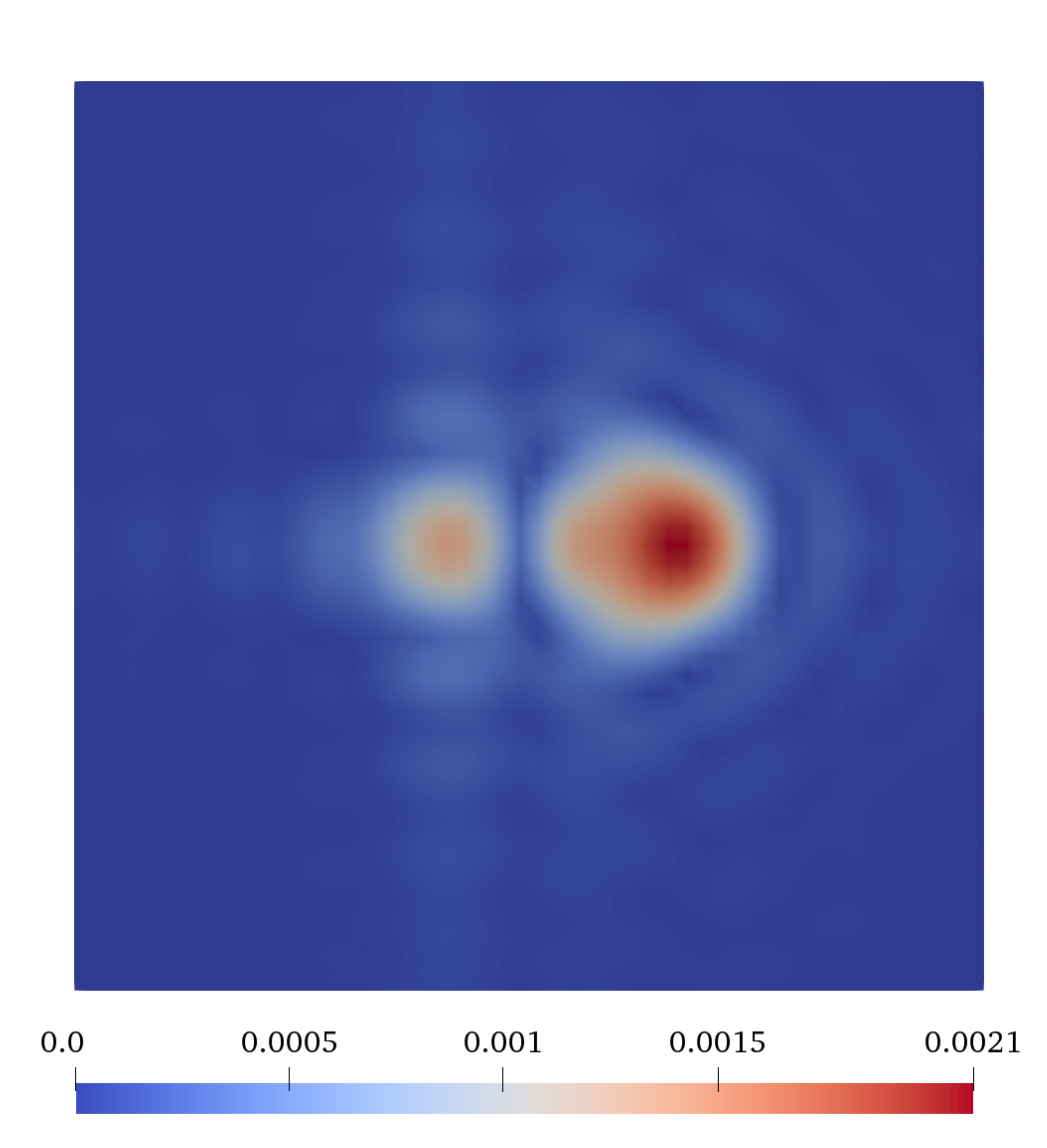}}
\end{center}
\caption{The absolute value of the simulated data with 15\% noise for the
case of two inclusions listed as number 4 in Table \protect\ref{tab1}. Here $%
k=16.2.$ a) The measured data $f_{noisy}\left( \mathbf{x},k\right) ,$ b) the
propagated data $g_{0,noisy}\left( \mathbf{x},k\right) ,$ c) the function $%
\protect\phi _{0,noisy}.$ Observe that the data propagation procedure helps
to separate these two inclusions. Also, it is clear from b), c) that the
left inclusion has a larger dielectric constant and a smaller size than the
right inclusion.}
\label{fig:2inc_data}
\end{figure}

Until now we have considered only the case of a single inclusion. The case
of two inclusions, which is listed as number 4 in Table \ref{tab1}, is very
similar. The absolute value of simulated data with noise $f_{noise}(\mathbf{x%
},k)$, the propagated data $g_{0,noisy}(\mathbf{x},k)$, and the function $%
\phi _{0,noisy}(\mathbf{x},k)$ for two inclusions and the wavenumber $k=16.2$
are displayed on Fig. \ref{fig:2inc_data}.

\begin{figure}[tbp]
\begin{center}
\subfloat[\label{fig:c_true2inc}]
{\includegraphics[width=0.5\textwidth]{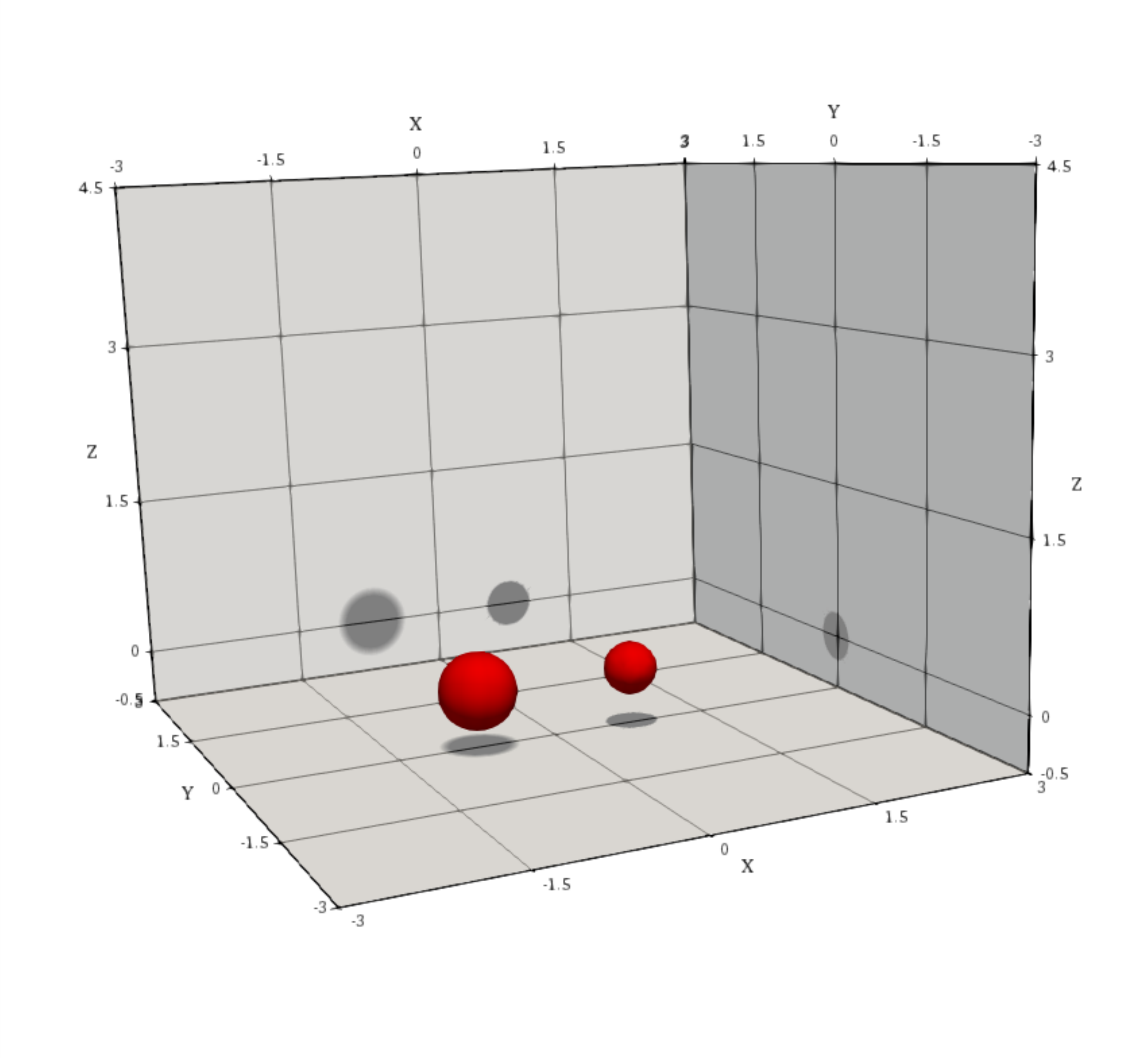}} 
\subfloat[%
\label{fig:c_comp2inc}]{\includegraphics[width=0.5\textwidth]{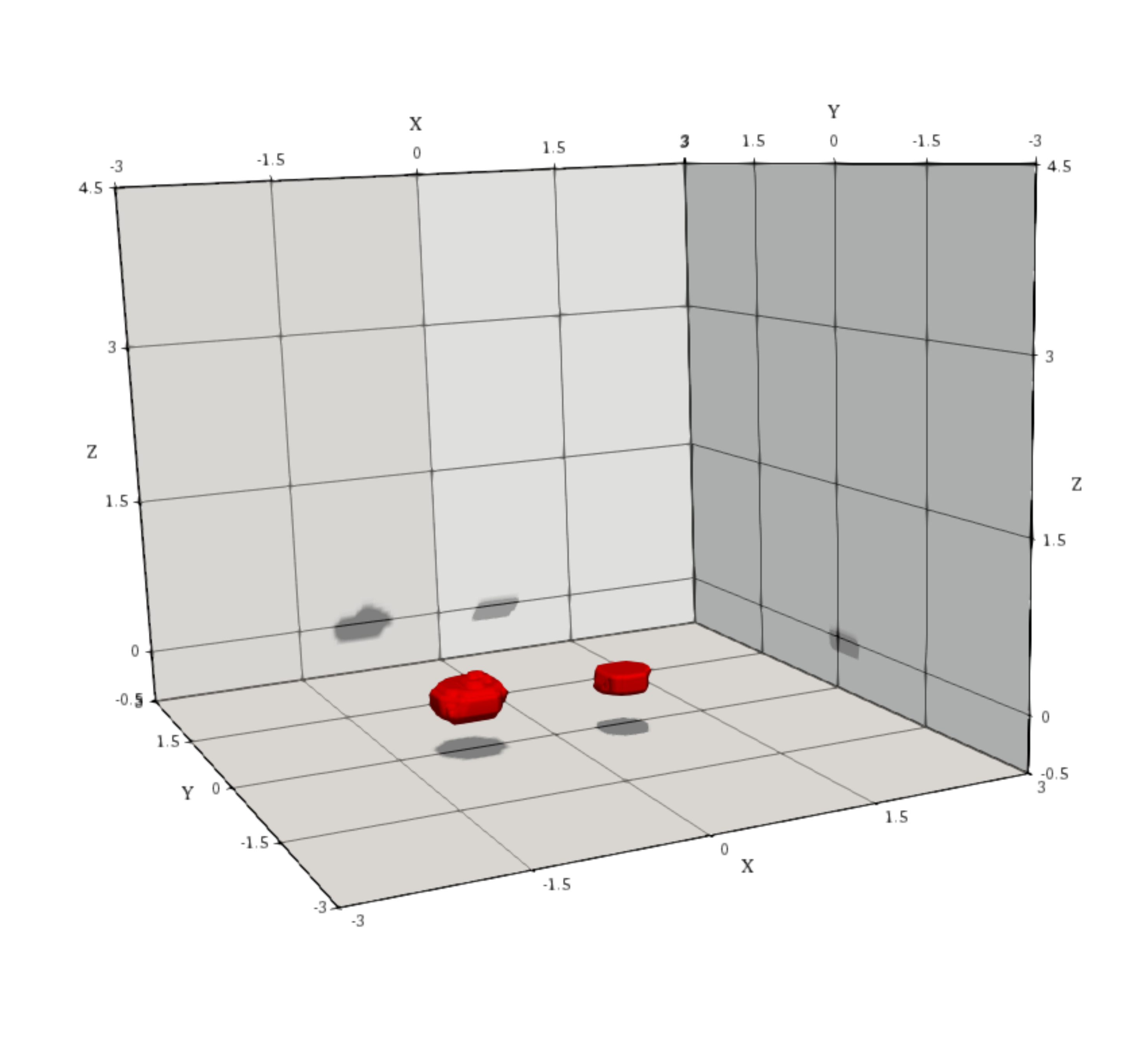}}
\end{center}
\caption{Reconstruction result for two inclusions: exact image (a) and
computed image (b)}
\label{fig:c_2inc}
\end{figure}

Looking at the original data of Fig. \ref{fig:f_2inc}, we cannot clearly
distinguish these two inclusions. However, Figures \ref{fig:g_2inc} and \ref%
{fig:q_2inc} show that these two inclusions can be clearly separated after
the data propagation procedure. Furthermore, these figures also indicate
that the left inclusion has a larger dielectric constant and a smaller size
than the right one, which is true. The reconstruction results of Fig. \ref%
{fig:c_2inc} reflect this fact too. Here, the locations of both inclusions
are computed accurately and the larger inclusion appears larger in the
reconstructed image \ref{fig:c_comp2inc}. The values of $c_{comp}$ in both
inclusions are also computed with a good accuracy, see Table \ref{tab3}.
This result is obtained using the same parameters as in the case with a
single inclusion.

\begin{table}[tbp]
\caption{Reconstruction results}
\label{tab3}
\begin{center}
\begin{tabular}{ccccc}
\hline
Inclusion number & Exact coef. $c_{exact}$ & Radius $r$ & Computed coef. $%
c_{comp}$, error & Location \\ \hline
1 & 3 & 0.3 & 3.17, 5.7\% & 0.01 \\ 
2 & 3 & 0.5 & 2.88, 4.0\% & 0.01 \\ 
3 & 5 & 0.3 & 5.15, 3.0\% & -0.09 \\ 
4.1 & 7 & 0.3 & 6.36, 9.0\% & 0.01 \\ 
4.2 & 3 & 0.5 & 2.99, 0.3\% & 0.01 \\ \hline
\end{tabular}%
\end{center}
\end{table}

\end{document}